\documentclass[12pt]{article}
\usepackage{amsfonts}
\usepackage{amssymb}
\usepackage{float}
\usepackage{graphicx}
\usepackage{caption}
\usepackage{subcaption}

\graphicspath{ {./newimages/} }
\graphicspath{ {./newimages/} }
\graphicspath{ {./newimages/} }

\input{tcilatex}
\begin{document}

\begin{center}
\textbf{The Kind of Silence: Managing a Reputation for Voluntary Disclosure
in Financial Markets}

\bigskip

Miles B. Gietzmann\footnote{%
Bocconi University; miles.gietzmann@unibocconi.it} and Adam J. Ostaszewski%
\footnote{%
London School of Economics; a.j.ostaszewski@lse.ac.uk; Corresponding author}

\bigskip
\end{center}

\textbf{Abstract. }In a continuous-time setting we investigate how the
management of a firm controls a dynamic choice between two generic voluntary
disclosure decision rules (strategies): one with full and transparent
disclosure termed \textit{candid}, the other, termed \textit{sparing}, under
which values only above a dynamic threshold are disclosed. We show how
management are rewarded with a reputational premium for being candid. The
candid strategy is, however, costly because the sparing alternative shields
the firm from potential downgrades following low value disclosures. We show
how parameters of the model such as news intensity, pay-for-performance and
time-to-mandatory-disclosure determine the optimal choice of candid versus
sparing strategies and the optimal times for management to switch between
the two. The private news updates received by management are modelled
following a Poisson arrival process, occurring between the fixed (known)
mandatory disclosure dates, such as fiscal years or quarters, with the news
received generated by a background Black-Scholes model of economic activity
and of its partial observation. The model presented develops a number of
insights, based on a very simple ordinary differential equation (ODE)
characterizing equilibrium in a piecewise-deterministic model, derivable
from the background Black-Scholes model. It is shown that in equilibrium
when news intensity is low a firm may employ a \textit{candid} disclosure
strategy throughout, but will otherwise switch (alternate) between periods
of being \textit{candid} and periods of being \textit{sparing} with the
truth \textbf{(}or\ the\ other\ way\ about\textbf{)}; that is, we are able
to characterize when in equilibrium a candid firm will switch to adopting a
sparing strategy. The model illustrates how parameters such as time to
mandatory disclosure, news intensity and pay-for-performance may drive such
switching behaviour. In summary: \textit{with constant pay-for-performance
parameters, at most one switching can occur.}

\textbf{Mathematics Subject Classification}: Primary 91G50; 91-10; 91B38;
91B44;

Secondary 91G80; 93E11; 93E20; 93E35; 60G35; 60G25

\textbf{Keywords: }Asset-price dynamics; voluntary disclosure; dynamic
disclosure policy; Markov piecewise-deterministc modelling; corporate
transparency reputation.

\section{Introduction}

When investors value firms, they not only base their inferences upon what
news signals managers make public, but also on the likelihood that
management may be hiding other news. A legal environment with high penalties
reduces the chances managers will hide very bad news signals, but a constant
concern for investors is: do management release early warning signals of
potentially less severe bad news in a timely fashion, or do they hide it in
an underhand way in the hope conditions will evolve differently and only
disclose when potential legal liability arises? (Marinovic and Varas [MarV].)

The Dye [Dye] model of voluntary disclosure addressed this issue in 1985 in a
static setting and derives equilibrium conditions under which management
adopts a sparing approach with a threshold-disclosure strategy, when
deciding whether or not to voluntarily disclose new information ahead of a
mandatory disclosure date. Institutional features, such as news-arrival
rates and time to the mandatory disclosure date, cannot be modelled in such
a static setting. In response Beyer and Dye [BeyD] (in 2012) develop a two-period
model in which managers may make a voluntary disclosure in order to build a
reputation for being \textit{candid} -- forthright (or `forthcoming') -- by
always faithfully disclosing their updated information. By contrast managers
could exploit their asymmetric endowment of information and\textbf{\ }only
make voluntary disclosures if their received signal is high enough, as
defined by a valuation threshold (cutoff). Such behaviour will here be
termed \textit{sparing}\footnote{%
Candid and sparing reporting strategies are termed \textit{forthcoming} and
\textit{strategic }by Beyer and Dye.}. For the two-period model Beyer and Dye show why managers' concerns in the
first period -- for how investors form second-period inferences (based on
observed first-period dividend outcomes) -- affect their voluntary
disclosure strategy. Their model predicts a diversity in management
strategies, as in equilibrium some managers will choose a strategy which
leads investors to assign a high probability that they will behave
sparingly, while others choose to be candid. An insight from this model is
that a voluntary disclosure strategy may be used to influence future
firm-value over and above the direct effect of disclosure of any
idiosyncratic signals of value. That is, establishing a reputation for being
candid\textbf{\ }at\ times\textbf{\ }shifts firm-value upwards, over and
above the direct discounted present-value of the most recent signals of
value, precisely because investors now assign a reduced likelihood for
management hiding bad news.

The structure of the paper is as follows. We discuss related literature in Subsection 1.1.
Section 2 presents our main
findings. Section 3 develops the theory of the optimal sparing-disclosure
threshold in a continuous framework, for which the main optimization tool
comes from control theory and relies on the Pontryagin Principle. Examples
are given in Section 4. Proofs of theorems are in Section 5 with some
technical results relegated to Appendix A and further, more routine, calculations to Appendix B, which ends with a \textit{symbols list}. Section 6 briefly considers multi-switching, although our main interest is in single-switching. We present conclusions in Section 7. We note that all valuations are viewed as {\textit{discounted}}.

\subsection{Related Literature on the Disclosure Dynamics}

Earlier models of dynamic disclosure have taken various paths, as follows.
Beyer and Dye comment that cheap talk models are concerned with reputational
formation rather than establishing a reputation for timely disclosure. Also
disclosure has the feature of a variable that must be binary and in practice
managers can disclose untruthfully, which is excluded in the Dye disclosure
model (and by us below). Another path, followed by Acharya et al. in
[AchMK], has been to investigate clustering of corporate disclosures driven
by market news. In that paper management learns only one piece of
information and must decide when to disclose it if other market news events
are occurring. Guttman et al. [GutKS] extend the Dye model to three periods
and demonstrate a differential equilibrium response to end of first or
second period disclosures. Marinovic and Varas [MarV] develop a
continuous-time disclosure model, but one in which there are exogenous costs
of disclosure. They show how litigation risk affects disclosure. They
explain how litigation stimulates transparency in financial markets and,
furthermore, how bad news can crowd out good news, leading to investors
being less skeptical about corporate disclosures. Bertomeu, Marinovic, Terry
and Varas [BerMTV] develop a model of dynamic concealment which is in part
motivated by a desire to develop a structural estimation model applied to
publicly observable annual earnings forecasts. In their dynamic disclosure
model, non-disclosure corresponds to firms not issuing an earnings forecast
(76 p.c. of cases) and so their model can be viewed as a between-fiscal
periods model of disclosure. In contrast, we focus on how investors should
rationally update value estimates given non-disclosure in the period between
fiscal dates (years or quarters). Since the institutional setting is
different for between- versus within-fiscal period settings, we provide a
\textit{disclosure market-microstructure}, which complements their
multi-year modelling. In our intra-period (within-period) model end-effects
arise in disclosure behaviour which may not be present in an inter-period
(between-periods) analysis, if there is no modelling of firm exit or
failure. Einhorn and Ziv [EinZ] develop a dynamic model, but one in which
disclosure generates a cost, and unlike the Dye model does not naturally
allow equilibrium risk-neutral pricing. The model closest in spirit to ours
is the Beyer and Dye [BeyD] model in which managers choose between a
strategy of either being \textit{candid} (always disclosing privately
observed news signals) or alternatively \textit{sparing}, only disclosing when it is in their own self-interest. We extend their two-period model to a
continuous-time setting and model the possibility of reputationally based
switching between the two reporting behaviours, driven by news-flow,
pay-for-performance, and or proximity to mandatory disclosure dates. Our
model shows that a managerial strategy of being candid over some
time-interval need not be a time-invariant strategy, instead it should only
be interpreted as a period-specific strategy. That is, our model calls into
question type-based models that assume managers have some deep intrinsic
desire to be of reputable type in all periods; instead, our model shows how
in equilibrium, rational, previously reputable, managers may be prepared to
burn reputation (by diverging from initial expected behavior) in a later
period. This suggests that it is rational for investors to condition beliefs
not just on history, but also on forward-looking variables such as remaining
tenure.

\section{Model generalities and findings}

Consider a firm whose financial state $X_{t}$ evolves in continuous time
according to a Black-Scholes model with periodic mandatory disclosure dates
and with interim intermittent capability of voluntary disclosures of the
next mandatory expected financial report. This would be based on partial
observation $Y_{t}$ of the financial state $X_{t}$.

Consider two possible reporting behaviours executed by the firm management
at any time $t$ when $Y_{t}$ is observed:

candid (faithful) reporting -- reporting the observed value, as seen, i.e.
unconditionally;

sparing (threshold) reporting -- reporting only the value observed when
above a time-$t$ dependent threshold, i.e. conditionally.

These behaviours are both capable of being applied at any one time, i.e.
leading to their use in some combination, and are assumed to be both
truthful and prompt (i.e. without delay). We also assume that managers
cannot \textit{credibly} assert absence of information arriving at time $t$
(i.e. absence of knowledge of $Y_{t}$). Furthermore, no\ evidence of an
undisclosed observation is retained. We admit no further sources of
information about $Y_{t};$ modelling with the inclusion of further sources
is touched on in [GieOS].

\textit{Sparing} here is used in the sense of being economical with
information delivery as in `economical with the truth' or `actualit\'{e}'%
\footnote{%
The phrase `economical with the truth', though it dates back to 1897, was
not common parlance in the UK till Robert Armstrong's reference to it -- in
defence of his stance during the Spycatcher trial in Australia in 1986,
resurfacing in 1992 in the Arms-for-Iraq case.}, sometimes called strategic.
It is of course a \textit{foundational question} whether it is suboptimal to
withold information. An early finding in the disclosure literature, provided
by Grossman and Hart [GroH] (in 1980) and Grossman [Gro] (in 1981), has become known as the
\textit{unravelling result}. It suggested that withholding information would
lead investors to discount the valuations, thus incentivising a firm to make
a full disclosure in order to restore the value.

The contribution of Dye [Dye] (in 1985) was to provide, in a discrete framework with
one interim date (say at some time $s$ between two mandatory disclosure
dates of $0$ and $1$), a rationale for why this `full disclosure
unravelling' result might not occur at the interim date $s,$ and to supply
an equation uniquely determining the resulting market discount in value in
an equilibrium framework. The market discount is an appropriate \textit{%
weighted average} that combines the possibility that management lacks fresh
information with the possibility that management may hide information which
if disclosed would have led to an even larger discount (i.e. below the
weighted average).

Dye's paradigm for valuing a non-disclosing firm may be characterized by an
amended statement of the Grossman-Hart paradigm as follows

\bigskip

\noindent \textbf{Minimum Principle} ([OstG], cf.
[AchMK] -- their Prop. 1). \textit{In equilibrium the market
values the firm at the least level consistent with the beliefs and
information available to the market as to its productive capability.}

\bigskip

See\textbf{\ }also Section 3. This result carries some detailed implications
to which we return later. But the principle already suggests intuitively
that if the management reporting behaviour is believed by the market to be
at times candid, then in equilibrium the weighted average valuation may at
times move upwards, by placing less weight on the chance of poor observation
being witheld.

We will demostrate the validity of such a suggestion in the continuous-time
context of the firm as described above, by creating a continuous analogue of
Dye's argument in which the Poisson \textit{arrival rate} of the observation
time of $Y_{t}$ is $\lambda $ and assuming management can report in a
sparing mode (relative to an optimally generated threshold) with a
probability $\pi _{t}$ at time $t$ when simultaneously the market believes
(in equilibrium) that the selected probability is indeed $\pi _{t}.$

Management choice of $\pi _{t}$ is motivated through the maximization of an
appropriate objective function rewarding in proportion to a factor $\alpha
_{t}$ the instantaneous firm value and penalizing in proportion to a factor $%
\beta _{t}$ the instantaneous \textit{value-differential} (value relative to
that derived from sparing behaviour executed throughout all time). The
parameter $\kappa _{t}=1-\alpha _{t}/\beta _{t}$ emerges as significant (see
Section 3.3\textbf{).}

The aim of the penalty term is to provide a tension between sparing and
candid reporting: adopting candour throughout would require more\ frequent\
disclosure of potentially bad news, which could result in larger falls in
firm value, i.e. over and above falls that resulted from continued sparing
silence (non disclosure).

We discuss our findings in this section, leaving details of the optimization
and proofs to the next and later section. Our first surprising finding is
that an optimal disclosure behaviour is of \textit{bang-bang }type which
will always \textit{switch }(alternate) between intervals of constancy with
only $\pi =0$ or $\pi =1,$ i.e. a mixed strategy is ruled out in the
following theorem. (See Appendix A for the stronger statement in Theorem
1S.)

\bigskip

\noindent \textbf{Theorem 1 (Non-mixing Theorem). }\textit{When }$\alpha
_{t},\beta _{t}$ \textit{\ are constant:}\textbf{\ }\newline
\textit{A mixing control with }$\pi _{t}\in (0,1)$\textit{\ is non-optimal
over any interval of time.}

\textit{\bigskip }

The theorem agrees with empirical findings (due to Grubb [Gru]) that after
an announcement management is observed to follow initially either candid or
sparing behaviour but not a mixture.

Accordingly, we study $\pi _{t}\in \{0,1\}.$ In particular, we study the
possible occurrence of an \textit{initially candid (candid-first) equilibrium%
} in which management at first adopt candid behaviour out of which they
switch after some time $\theta ,$ the \textit{switching time}, in favour of
a sparing policy, and also \textit{initially sparing} (\textit{%
sparing-first) equilibrium} in which management at first adopt sparing
behaviour out of which after some time they switch in favour of candid
behaviour. This provides a model framework for \textit{empirical detection
of regime change} (disclosure policy change): cf. L\o kka [Lok]. Theorems
2a and 2b identify both the \textit{location} of the switching time of such
an equilibrium policy and the attendant necessary and sufficient \textit{%
existence conditions} guaranteeing an equilibrium. These theorems are
followed by a clarifying discussion concerning the location conditions.

We stress that having merely a characterization of the location condition is
not adequate. The technical nature of the existence conditions emerges from
a Hamiltonian analysis (Section 3.4 below) in which the Pontryagin Principle
relies on Theorem 1 (the non-mixing) in supplying a necessary and sufficient
optimality condition.

Our findings refer to a decreasing discount function $h(t)$ (responsible for
the rate of fall in values when continued absence of disclosure is
attributed to sparing behaviour -- see Section 3.2\ equation (cont-eq)) and
to its integral%
\[
g(t)=\int_{0}^{t}h(u)\text{ }\mathrm{d}u.
\]

\noindent \textbf{Theorem 2a (Single switch equilibrium location and}
\textbf{existence for an initially candid strategy}).\newline
\textit{Assume that }$\alpha _{t},\beta _{t}$\textit{\ are constant and }$%
0<\kappa <1$\textit{\ for }$\kappa :=1-\alpha _{t}/\beta _{t}$\textit{.}%
\newline
\textit{In an equilibrium, if such exists, in which }$\pi =0$ \textit{on} $%
[0,\theta )=0$\textit{\ and }$\pi =1$\textit{\ on }$[\theta ,1],$\ \textit{%
the uniquely optimal switching time }$\theta $\textit{\ solves}%
\[
\kappa _{1}:=\kappa e^{\lambda g(\theta )}=1,\text{ i.e.}\qquad g(\theta
)=-\left. \log \kappa \right/ \lambda .
\]%
\textit{For given }$\lambda $ \textit{this equation is solvable for large
enough }$\kappa ,$ \textit{in fact iff}%
\begin{equation}
1>\kappa \geq e^{-\lambda g(1)}.  \tag{cand}
\end{equation}%
\textit{Such an equilibrium exists iff the unique switching time }$\theta $
\textit{satisfies}%
\[
g(\theta )/(\theta h(0))>\frac{\log \kappa ^{-1}}{\kappa ^{-1}-1}.
\]%
\textit{In such an equilibrium, the unique switching time }$\theta $ \textit{%
satisfies }%
\[
\theta ^{\prime }(\lambda )=\frac{\log \kappa }{\lambda ^{2}h(\theta
(\lambda ))}<0.
\]%
\textit{So larger news-arrival rates }$\lambda $\textit{\ create shorter
periods of initial candid behaviour.}

\bigskip

\noindent \textbf{Remark.} The left-hand side term of the existence
condition above is monotonically decreasing from $1$ down to $g(1)/h(0)$, as
in the red graph in Figure 1 below; its lowest value is dictated by $\sigma
, $ a volatility measure. The right-hand side ranges monotonically from $0$
to $1;$ thus a fixed value of $\kappa $ supplies a {value to the right-hand
side and is illustrated in green below for a choice which allows all values $%
\theta $ to satisfy the inequality here (with other choices restricting the $%
\theta $ range).}
\begin{figure}[tbp]
\begin{center}
\includegraphics[height=3.5cm]{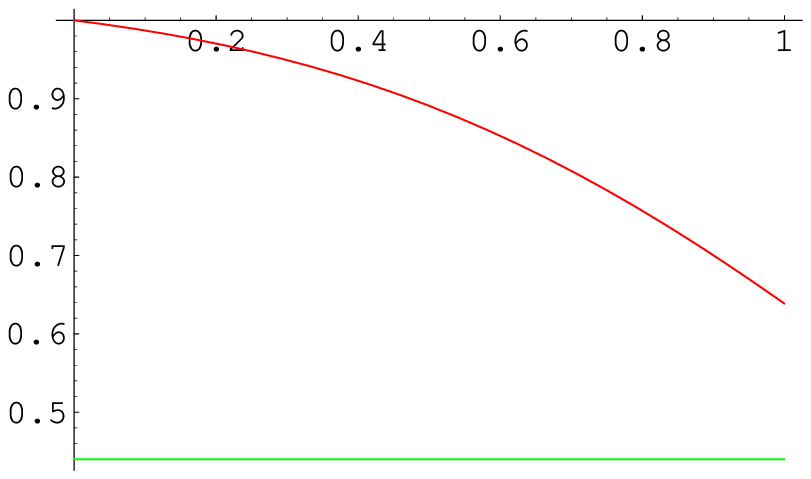}
\par
{Figure 1: Pay-for-performance bounds}
\end{center}
\end{figure}

Matters are more complicated in an equilibrium that is initially sparing.

\bigskip

\noindent \textbf{Theorem 2b (Single switch equilibrium location and}
\textbf{existence for an initially sparing strategy}).\newline
\textit{Assume that }$\alpha _{t},\beta _{t}$\textit{\ are constant and }$%
0<\kappa <1$\textit{\ for }$\kappa :=1-\alpha _{t}/\beta _{t}$\textit{.}%
\newline
\textit{In an equilibrium, if such exists, in which }$\pi =1$ \textit{on} $%
[0,\theta )=0$\textit{\ and }$\pi =0$\textit{\ on }$[\theta ,1],$\ \textit{%
the uniquely optimal switching time }$\theta $\textit{\ solves the
first-order condition}%
\[
1=\theta +\frac{\kappa ^{-1}-1}{\lambda h(\theta )},\text{ equivalently }%
(1-\theta )h(\theta )=\frac{\kappa ^{-1}-1}{\lambda }.
\]%
\textit{For given }$\lambda ,$ \textit{this equation is solvable for large
enough }$\kappa ,$ \textit{in fact iff}%
\begin{equation}
1>\kappa \geq \lbrack 1+\lambda h(0)]^{-1}.  \tag{spar}
\end{equation}

\textit{Such an equilibrium exists iff the unique switching time }$\theta $
\textit{satisfies}%
\begin{eqnarray*}
(1-\theta )h(\theta )(-[\log h(\theta )/h(0)]/g(\theta )) &>&\kappa ^{-1}-1,%
\text{ (a lower bound for }\theta \text{),} \\
-[\log h(\theta )/h(0)]/g(\theta ) &>&\lambda \text{ (a bound on }\lambda
\text{ in terms of }\theta \text{).}
\end{eqnarray*}%
\textit{In such an equilibrium, the unique switching time }$\theta $\textit{%
\ satisfies }%
\[
\theta ^{\prime }(\lambda )=\frac{(1-\theta )h(\theta )}{\lambda \lbrack
h(\theta )-(1-\theta )h^{\prime }(\theta )]}>0.
\]%
\textit{So larger news-arrival rates }$\lambda $\textit{\ create longer
periods of initial sparing behaviour.}

\bigskip

The theorems expose a fact of direct relevance to empirical study: that with
the same parameter values it may happen that both equilibrium types coexist,
as in Figures 2a and 3. (Conditions on parameter values $\kappa ,\lambda $\
permitting this can be derived numerically from the conditions of Theorems
2a and 2b.) In particular, the two conditions (cand) and (spar) on $\kappa $
with fixed $\lambda ,$ may both hold simultaneously: indeed, since the map $%
\lambda \mapsto $ $e^{\lambda g(1)}$ is convex, there is a unique $\lambda
=\lambda _{\text{crit }}>0$ such that%
\[
e^{\lambda g(1)}=1+\lambda h(0),
\]%
and%
\[
e^{-\lambda g(1)}\lessgtr (1+\lambda h(0))^{-1}\text{ according as }\lambda
_{\text{crit}}\lessgtr \lambda .
\]%
Thus, for $\lambda <\lambda _{\text{crit}},$ both conditions are met when $%
\kappa $ satisfies (cand):
\[
(1+\lambda h(0))^{-1}<e^{-\lambda g(1)}<\kappa .
\]
In contrast to the coexistence of equilibria in which switching occurs,
there does exist a constantly candid equilibrium (i.e. with no switching):

\bigskip

\noindent \textbf{Qualitative Corollary.} \textit{For small enough }$\lambda
,$ \textit{candid (unconditional) disclosure throughout the period of
silence is an equilibrium policy.}

\bigskip

For proof: see Corollary 2 in Section 5. The location of optimal switches,
assuming they correspond to an equilibrium (requiring additional
conditions), can be pursued in generality. We identify the consequent
generalization as this leads to yet another surprising finding stated in the
Corollary below. For the proof\ of Theorems 3a and 3b see Appendix B.

\bigskip

\noindent \textbf{Theorem 3a (General switching equation). }\textit{Assume
that }$\alpha _{t},\beta _{t}$\textit{\ are constant and }$0<\kappa <1$%
\textit{\ for }$\kappa :=1-\alpha _{t}/\beta _{t}$\textit{.}\textbf{\ }%
\textit{Assume further that the }$n$\textit{\ consecutive switches located
at times }$0=\theta _{0}<\theta _{1}<\theta _{2}<...<\theta _{n}<1=\theta
_{n+1}$ \textit{are selected optimally. Then, for }$1\leq i\leq n,$\textit{\
}%
\begin{eqnarray*}
g\left( \theta _{i}+\frac{1-\bar{\kappa}_{i-1}e^{\lambda g(\theta _{i-1})}}{%
\bar{\kappa}_{i-1}e^{\lambda g(\theta _{i-1})}\lambda h(\theta _{i})}\right)
&=&g(\theta _{i})-(\log \bar{\kappa}_{i-1})/\lambda -g(\theta _{i-1}),\text{
if }\pi =1\text{ on }[\theta _{i-1},\theta _{i}), \\
e^{-\lambda g(\theta _{i})} &=&\bar{\kappa}_{i},\text{ if }\pi =0\text{ on }%
[\theta _{i-1},\theta _{i}),
\end{eqnarray*}%
\textit{where the parameters }$\bar{\kappa}_{j}$\textit{\ defined below
depend on }$\{\theta _{k}:k\leq j-1\}.$%
\[
\bar{\kappa}_{j}:=\kappa \gamma _{j}\qquad (\text{for }\gamma _{j}:=\gamma
_{\theta _{j}}).
\]%
\bigskip

\noindent \textbf{Theorem 3b (Multiple switching locations). }\textit{Assume
that }$\alpha _{t},\beta _{t}$\textit{\ are constant and }$0<\kappa <1$%
\textit{\ for }$\kappa :=1-\alpha _{t}/\beta _{t}$\textit{.} \textit{The
sequence of solutions to the switching equation defines the switching times
according to the recurrence}%
\begin{eqnarray*}
\theta _{i+1} &=&\theta _{i}+\frac{\bar{\kappa}_{i-1}^{-1}e^{-\lambda
g(\theta _{i-1})}-1}{\lambda h(\theta _{i})},\text{ if }\pi =1\text{ on }%
[\theta _{i-1},\theta _{i}), \\
g(\theta _{i}) &=&\log \bar{\kappa}_{i}^{-1}/\lambda ,\text{ if }\pi =0\text{
on }[\theta _{i-1},\theta _{i}),
\end{eqnarray*}%
\textit{unless }$i=n,$\textit{\ so that }$\theta _{i+1}=1,$\textit{\ in
which case }$\theta _{i}$\textit{\ is determined by the equation}%
\[
\theta _{i+1}=1=\theta _{i}+\frac{\bar{\kappa}_{i-1}^{-1}e^{-\lambda
g(\theta _{i-1})}-1}{\lambda h(\theta _{i})}.
\]%
\textit{If }$\theta _{1}$\textit{\ is a right endpoint of an interval where }%
$\pi =0$\textit{, then}%
\[
g(\theta _{1})=-\log \kappa /\lambda .
\]%
\textit{Furthermore, the sequence }$\bar{\kappa}_{i}$\textit{\ is (weakly)
decreasing with alternate members strictly decreasing.}

\bigskip

\noindent \textbf{Corollary (Candid-first single switching).}
\textit{Assume that }$\alpha _{t},\beta _{t}$\textit{\ are constant and }$%
0<\kappa <1$\textit{\ for }$\kappa :=1-\alpha _{t}/\beta _{t}$\textit{.}%
\newline
\textit{If }$\pi =0$\textit{\ on }$[0,\theta _{1}),$\textit{\ so that }$%
\gamma _{1}=\gamma _{0}=1,$\textit{\ then }$\pi =1$\textit{\ on }$[\theta
_{1},\theta _{2})$\textit{\ and so for }$i=2,$\textit{\ as }$\bar{\kappa}%
_{1}=\kappa =e^{-\lambda g(\theta _{1})},$\textit{\ }%
\[
\theta _{3}=\theta _{2}+\frac{1-\bar{\kappa}_{1}e^{\lambda g(\theta _{1})}}{%
\bar{\kappa}_{1}e^{\lambda g(\theta _{1})}\lambda h(\theta _{2})}=\theta
_{2},
\]%
\textit{a contradiction to }$\theta _{2}<\theta _{3}$\textit{. Consequently,
there cannot be a further switching from sparing to candid mode.}

\bigskip

A similar result appears to be supported by numerical\ analysis for a
sparing-first equilibrium policy, albeit Theorem 3b (on its own, i.e.
without invoking equilibrium conditions) implies by a similar argument that
if\textit{\ }$\pi =1$\textit{\ }on\textit{\ }$[0,\theta _{1}),$ then $\theta
_{4}=\theta _{3},$ i.e. at most two switchings can occur.

In summary, this section has characterized how tractable single-switching
conditions can be derived. The issue of equilibrium selection (which of the
sparing-first or the candid-first) must rest on the underlying assumption
that the market has found its way to one or other of the two by some
evolutionary game-theoretic mechanism; for a standard textbook view of the
latter, see e.g. Weibull [Wei].

\textbf{\bigskip }

\noindent \textbf{Remark. }Optimal multi-switching becomes possible when a
time-varying $\kappa _{t}$\ replaces the constant $\kappa ;$\ this is
particularly easy to arrange in the case of a piecewise-constant $\kappa
_{t} $\ with constancy on each inter-switching interval (interval between
successive switching points): see Example 3 in Section 3.

We have seen above that single-switching should be regarded as natural in
the constant $\kappa $\ context and not just a stylized model choice.
Moreover, single-switching provides the pragmatic, preferred equilibrium
choice by an appeal to \textit{focal-point} (Schelling-point) equilibrium
selection -- for a standard textbook view of which, see e.g. [FudT].

\section{The sparing disclosure threshold}

For $X_{s}$ the financial state and $Y_{s}$ its observation assume the
regression function $m_{s}(y):=\mathbb{E}_{t=0}[X_{1}|Y_{s}=y]$ is
increasing. Then the opimal threshold $\gamma _{s}$ is uniquely determined
and has three properties, the first of which in (i) below implies the
Minimum Principle of Section 2.

(i) \textbf{Minimum Principle: }[OstG]. The
valuation function
\[
W(\gamma ):=\mathbb{E}_{0}[X_{1}|ND_{s}(\gamma )]
\]%
has a unique \textit{minimum} at $\gamma =\gamma _{s};$

(ii) \textbf{Risk-neutral Consistency Property: }$\gamma _{s}$ is the unique
value $\gamma $ such that%
\[
\mathbb{E}_{0}[X_{1}|Y_{0}]=\tau _{\text{D}}\mathbb{E}_{0}[X_{1}|Y_{t}\geq
\gamma \text{ disclosed}]+\tau _{\text{ND}}m_{s}(\gamma ),
\]%
with $\tau _{\text{D}}$ is the (time $t$) probability of disclosure occuring
at $s.$ This is highly significant, in that the valuation at time $0$
anticipates the potential effects of a voluntary disclosure at the future
interim date $s.$ In brief, the approach is consistent with the principles
of \textit{risk-neutral valuation}; for background see Bingham and Kiesel [BinK], Chap. 6.
In particular, the risk-neutral valuation is a {\textit{martingale}},
constructed via iterated expectations from $\gamma_s:=\mathbb{E}_s[X_1|Y_1]=\mathbb{E}_s[\gamma_1]$
-- see [GieO] \S 2.3.

(iii) \textbf{Interim discount }From the perspective of time $0,$ in a model
with only three dates: $0,s,1,$ the Dye equation at time $s$ may be written:

\[
\frac{1-q_{s}}{q_{s}}(\gamma _{0}-\gamma _{s})=\int_{0}^{\gamma _{s}}F_{0}(u)%
\text{ }\mathrm{d}u=\int_{0}^{\gamma _{s}}(\gamma _{s}-u)^{+}\text{ }\mathrm{%
d}F_{0}(u),
\]%
and has the interpretation of a\textit{\ protective put-option} with strike $%
g$ against a fall in value at $t.$ Here $q_{s}$ is the probability that $%
Y_{s}$ is observed, and $F_{0}(u)=\mathbb{Q}(X_{1}\leq u|Y_{0}).$

The argument\ leading to these results is also sketched in the next section.

\subsection{Derivation of the Dye threshold equation}

Relocating the dates to $t<s<1,$ the interim discount $\gamma _{s}$, which
is also the threshold for announcements in equilibrium at time $s,$ is the
value $\gamma =\gamma _{s}$ which satisfies%
\begin{equation}
\gamma =\mathbb{E}^{\mathbb{Q}}(X_{s}|ND_{s}(\gamma ),\mathcal{F}_{t}),
\tag{indif}
\end{equation}%
with $\mathbb{Q}$ $=$ market's probability measure for all relevant events,%
\newline
$\mathcal{F}_{t}=$ market information at time $t$;\newline
$ND_{s}(\gamma )=$ event at time $s$ that no disclosure occurs when the
information is below $\gamma $;\newline
$RHS=$ market's expectation of value conditional no-disclosure $%
ND_{s}(\gamma _{s})$.

\bigskip

If management observe a value $\gamma _{s}$ at time $s,$ then they are
indifferent between disclosing the valuation as $\gamma _{s}$ and
withholding said information.

If the probability of information reaching management at time $s$ is $%
q_{s}=q,$ assumed exogenous and independent of the state of the firm, then
for $p:=1-q$
\begin{equation}
\gamma _{s}=\frac{p\gamma _{t}+q\int_{0}^{\gamma _{s}}xd\mathbb{Q}_{t}(x)}{p+%
\mathbb{Q}_{t}(\gamma _{s})},  \tag{cond}
\end{equation}%
as $\gamma_t=\mathbb{E}_t[\gamma_1]$, where the subscript indicates conditioning on $\mathcal{F}_{t}.$
Equivalently, we have%
\begin{equation}
p(\gamma _{t}-\gamma _{s}) =q\int_{0}^{\gamma _{s}}(\gamma _{s}-x)d\mathbb{%
Q}_{t}(x), \tag{put}
\end{equation}
\[
\mathbb{Q}_{t}((\gamma _{s}-X_{s})^{+}) =\mathit{protective\ put.}
\]
Here we may routinely evaluate this put using the Black-Scholes formula.

As above, rearrangement will show incorporation of future information:%
\begin{equation}
\gamma _{t}=\tau _{D}^{t}\cdot \mathbb{E[}x|D_{s}(\gamma _{s}),\mathcal{F}%
_{t}]+(1-\tau _{D}^{t})\gamma _{s},  \tag{r-n val}
\end{equation}%
as in risk-neutral valuation, where\newline
$D_{s}(\gamma )=$ event of time $s$ when values above $\gamma $ disclosed%
\newline
$\tau _{D}^{t}=$ market's evaluation at time $t$ of disclosure probability
at time $s.$

The presumption this far precluded the use of a candid strategy. If
management restricts application of the sparing (threshold-generated)
strategy to act with probability $\pi $ and the market likewise believes (in
equilibrium) that this probability is $\pi ,$ then in a period of silence:%
\begin{equation}
\gamma _{s}=\frac{p\gamma _{t}+\pi q\int_{0}^{\gamma _{s}}xd\mathbb{Q}_{t}(x)%
}{p+\pi \mathbb{Q}_{t}(\gamma _{s})}.  \tag{cond-$\pi $}
\end{equation}%
For $\pi =1$ (sparing) this reduces to the Dye equation.\newline
For $\pi =0$ (candid) this acknowledges that $\gamma _{s}=\gamma _{t}.$

\subsection{Equilibrium condition: continuous-time analogue}

We embed the three dates $t<s<1$ of the Dye model into the unit interval to
provide a continuous-time framework in which any future date $s>t$ can be
interpreted as a time at which the management have the opportunity to
disclose a forecast of value to the market. As in the Dye model, key here is
the creation of an ambiguity at time $s$, so that the market knows that
absence of a disclosure is caused either by absence of fresh endowment of
private information or by a management decision to withhold the private
information arriving at moment $s.$ With this aim we introduce a Poisson
process with intensity $\lambda $ whose jump at time $s,$ when privately
observed by management, determines that an observation of $Y_{s}$ occurs.
The market does not observe the jumps. Thus every moment now takes on the
character of an interim disclosure date and, depending on the disclosure
policy believed by the market to be implemented by management, absence of a
disclosure can mean no new observation or a withheld observation.

With the Poisson process in place, for $t<s$ take $q=q_{ts}=\mathbb{\lambda }%
(s-t)+o(s-t),$ employing the Landau little-o notation. Passage to the limit
as $s\searrow t$ yields:%
\begin{eqnarray*}
(1-q_{ts})(\gamma _{t}-\gamma _{s})+o(s-t) &=&q_{ts}\int_{z_{s}\leq \gamma
_{s}}(z_{s}-\gamma _{s})\text{ }\mathrm{d}\mathbb{Q}_{t}(z_{s}) \\
&=&-\mathbb{\lambda }(s-t)\int_{z_{s}\leq \gamma _{s}}\mathbb{Q}_{t}(z_{s})%
\text{ }\mathrm{d}z_{s}.\text{ (by parts)}
\end{eqnarray*}%
Dividing by $-\lambda (s-t)$:
\[
-(1-q_{ts})\frac{\gamma _{s}-\gamma _{t}}{(s-t)}=\lambda \int_{z_{s}\leq
\gamma _{s}}\mathbb{Q}_{t}(z_{s})\text{ }\mathrm{d}z_{s},
\]%
ignoring errors of order $o(s-t)/(s-t).$ With economic activity and the
noisy observation in a standard Black-Scholes setting, this yields%
\begin{equation}
-\gamma _{t}^{\prime }=\lambda \gamma _{t}h(t)\text{ with }h(t)=2\Phi
(\sigma (1-t)/2)-1,  \tag{cont-eq}
\end{equation}%
with $\sigma $ the \textit{(aggregate) volatility} (aggregating productive
and observation vols.); for the proof see [GieO].

This ODE\ is our continuous-time disclosure-equilibrium condition in any
period of silence (i.e. when the management make no disclosures). It
equilibrates in a period of silence between the market's ability to
downgrade the valuation below $\gamma _{t}$ and the management's potential
ability to upgrade the valuation were they to observe a greater value of $%
Y_{t}$ (cf. the weighted average discussed in Section 2). We refer to this
as the \textit{equilibrium ODE}.

Denoting successive public disclosures (voluntary or mandatory) generated
stochastically by $\tau _{0}=0<\tau _{1}<\tau _{2}<...<1$, and writing $N$
for their number so that $\tau _{N+1}=1,$ one has
\begin{eqnarray*}
\gamma _{t}^{\prime } &=&-\mathbb{\lambda }\gamma _{t}h(t-\tau _{i})\text{
for }\tau _{i}<t<\tau _{i+1} \\
\text{s.t. }\gamma _{\tau _{i}} &=&\text{disclosed value at time }t=\tau
_{i},
\end{eqnarray*}%
with $h(t-\tau _{i})$ the (per-unit of value, $\gamma $) the \textit{%
firm-specific,} \textit{instantaneous protective put} at times $t$ for $\tau
_{i}<t<\tau _{i+1}.$

The market valuation of the firm $\gamma _{t}$ is thus a \textit{%
piecewise-deterministic Markov} processes in the sense of Davis [Dav84, Dav93].

\subsection{Probabilistic strategy optimization}

The governing equation of our continuous-time version of the Dye model, the
equilibrium ODE, is based on the assumption that the manager's objective is
to achieve the highest possible valuation at all times $t$ preceding the
subsequent mandatory disclosure date. However, if management follow the
conditional threshold rule with probability $\pi _{t}$ and otherwise
disclose the observation candidly with probability $1-\pi _{t},$ then, as in
\textbf{(}cond-$\pi $\textbf{), }for an equilibrium strategy $\pi $\ the
corresponding valuation $\gamma _{t}=\gamma _{t}^{\pi }$ satisfies\textbf{:}%
\begin{equation}
\gamma _{t}^{\prime }=-\mathbb{\lambda }\gamma _{t}\pi _{t}h(t)\text{ with }%
h(t)=2\Phi (\sigma (1-t)/2)-1,  \tag{cont-eq-$\pi $}
\end{equation}%
where $t=0$ corresponds to the last public disclosure (after a change of
origin here, mutatis mutandis). We rescale the valuation so that $\gamma
_{0}=1.$

Consistently with this last equation, we will employ the notation: $\gamma
_{t}^{1}$ for its solution when $\pi \equiv 1$ (sparing policy applied
throughout), so that%
\[
\gamma _{t}=\gamma _{t}^{1}\text{ denotes the solution of }\gamma
_{t}^{\prime }=-\mathbb{\lambda }\gamma _{t}h(t)\text{ with }\gamma _{0}=1.
\]%
With comparison against this solution in mind, the manager is now induced to
maximize an objective in selecting $\pi $ so as to yield%
\begin{equation}
\max_{\pi }\mathbb{E}\int_{0}^{1}(1-\pi _{t})(\alpha _{t}\gamma _{t}-\beta
_{t}[\gamma _{t}-\gamma _{t}^{1}])\text{ }\mathrm{d}t.  \tag{obj-1}
\end{equation}%
As before, $t=0$ denotes the most recent time of disclosure and unit time is
left to the mandatory disclosure (\textit{time to expiry)}.

This objective includes a \textit{penalty} proportional to $(\gamma
_{t}-\gamma _{t}^{1}).$ The \textit{amended unravelling principle} of
Section 2 implies that introduction in a market equilibrium of candour
(candid reporting) will cause the valuation $\gamma _{t}$ to exceed $\gamma
_{t}^{1}$ and the aim of the penalty is to motivate management into
protecting the value of the firm from potential\ falls\ in value if a candid
strategy is followed for too long (i.e. from excessive use of a candid
position).

In equivalent form, the objective may be rewritten as%
\begin{equation}
\max_{\pi }\mathbb{E}\int_{0}^{1}(1-\pi _{t})\beta _{t}[\gamma
_{t}^{1}-\kappa _{t}\gamma _{t}]\text{ }\mathrm{d}t\text{ for }\kappa
_{t}:=1-\alpha _{t}/\beta _{t}.  \tag{obj-2}
\end{equation}%
It is thus natural to demand that for some proper interval of time%
\[
\gamma _{t}^{1}>\kappa _{t}\gamma _{t},
\]%
so we make the \textit{blanket assumption}
\[
0<\alpha _{t}/\beta _{t}<1,\text{ i.e. }0<\kappa _{t}<1,
\]%
which enables discounting of $\gamma _{t}$ by $\kappa _{t}$ to a level below
$\gamma _{t}^{1}.$

\subsection{Hamiltonian analysis: Pontryagin Principle}

We approach the maximization problem via the \textit{Pontryagin Maximum
Principle}, PMP, for which see [BreP] (esp. Ch. 7 on
sufficiency conditions for PMP), or the more concise textbook sketches in
[Lib], [Sas], or [Tro]. It is also possible to
establish the results below by solving the \textit{Bellman equation} along
the lines of Davis [Dav93, p. 165], a matter we hope to return to
elsewhere.

In a period of silence, the valuation is deterministic and so we formulate
optimisation in Hamiltonian terms. We apply a standard Hamiltonian approach
from control theory to maximizing the objective of the preceding section by
treating $\gamma _{t}$ as a state variable and $\pi _{t}$ as a control
variable. Denoting the co-state variable by $\mu _{t},$ the Hamiltonian is
\[
\mathcal{H}(\gamma _{t},\pi _{t},\mu _{t})=[(1-\pi _{t})\beta _{t}(\gamma
_{t}^{1}-\kappa _{t}\gamma _{t})]+\mu _{t}[-\mathbb{\lambda }\gamma _{t}\pi
_{t}h(t)],
\]%
by construction linear in $\pi _{t}.$ So with $\mu _{t}$ \textit{continuous}
and \textit{piecewise smooth}:%
\[
\mu _{t}^{\prime }=-\frac{\partial \mathcal{H}}{\partial \gamma _{t}}=\beta
_{t}\kappa _{t}(1-\pi _{t})+\mu _{t}\pi _{t}\mathbb{\lambda }h(t),\text{
with }\mu _{1}=0,
\]%
where we follow the \textit{c\`{a}dl\`{a}g convention }that $\pi _{t}$
right-continuous with left limits and satisfies $0\leq \pi _{t}\leq 1$. Thus%
\[
\mu _{t}^{\prime }-\mu _{t}\pi _{t}\mathbb{\lambda }h(t)=\beta _{t}\kappa
_{t}(1-\pi _{t}).
\]

We now apply the Pontryagin Principle. Evidently, concentrating only on
terms involving $\pi _{t}$ below, the Hamiltonian%
\[
\mathcal{H}(\gamma _{t},\pi _{t},\mu _{t})=...-\pi _{t}\cdot \beta
_{t}[\gamma _{t}^{1}-\kappa _{t}\gamma _{t}+\mu _{t}\mathbb{\lambda }\gamma
_{t}h(t)/\beta _{t}]
\]%
is maximized by setting $\pi _{t}$ at 0 or 1 according as%
\[
\gamma _{t}[\kappa _{t}-\mathbb{\lambda }h(t)\mu _{t}/\beta _{t}]<\gamma
_{t}^{1}\text{ or}>\gamma _{t}^{1},\text{ resp.}
\]%
It emerges that $\mu _{t}\leq 0$ (see Appendix A, Proposition 3) consistent
with its being interpreted as a penalty term in $\mathcal{H}$, so%
\[
\lbrack \kappa _{t}-\mu _{t}\mathbb{\lambda }h(t)/\beta _{t}]>0.
\]%
This gives rise to an \textit{optimal switching curve} and associated
optimality rule:
\[
\gamma _{t}^{\ast }:=\frac{\gamma _{t}^{1}}{\kappa _{t}-\mathbb{\lambda }%
h(t)\mu _{t}/\beta _{t}}>0.
\]

\textbf{Proposition 1 (Optimality Rule). }\textit{A} \textit{necessary and
sufficient for }$\pi $\textit{\ to be optimal is given by} \textit{the rule}:%
\[
\pi _{t}=\left\{
\begin{array}{cc}
1,\text{ suppress }X_{t}\text{ unless }X_{t}\geq \gamma _{t} & \text{ if }%
\gamma _{t}\geq \gamma _{t}^{\ast }, \\
0,\text{ reveal }X_{t} & \text{if }\gamma _{t}<\gamma _{t}^{\ast }.%
\end{array}%
\right.
\]

\textbf{Proof. }By the Non-mixing theorem, $\pi _{t}$ can only take the
values $0$ and $1$ and so by the Pontryagin Principle the optimality
condition above is \textit{necessary and sufficient}: the strong form of
Theorem 1 (see Appendix A) asserts that if $\gamma _{t}=\gamma _{t}^{\ast }$
on an interval of time, then $\pi _{t}=1$ on that interval. \hfill $\square $

\bigskip

A corollary of the above form of $\gamma _{t}^{\ast }$ now follows.

\bigskip

\textbf{Qualitative Corollary.} \textit{A large enough valuation }$\gamma
_{t}$\textit{\ allows sparing reporting, a small enough valuation }$\gamma
_{t}$\textit{\ encourages candid revelation.}

\bigskip

\noindent \textbf{Remark. }Evidently, the value of $\pi _{t}$\ is not
instantly observable, so management may at any instance of bad news (however
defined) hide it and so deviate from their prescribed equilibrium strategy.
However, systematic deviation of this sort is statistically observable and
so deviation leads to loss of reputation, removing the very means by which
the firm maintains a higher valuation, which in turn hurts the deviating
agent. We therefore assume that managers hold themselves to their prescribed
equilibrium strategy. For further background on the \textit{Bayesian
persuasion} aspect here, see Kamenica and Gentzkow [KamG].

\section{Examples of equilibrium behaviour}
In this section we give three examples of different equilibrium behaviour in
the form of graphs which include the switching curve derived in the
preceding section. The role of the switching curve is to confirm, by
Proposition 1, the optimality of the equilibrium valuation.
\subsection{\textbf{Example 1}}
Here $\pi =0$ initially (candid). Figure 2a below and Figure 3 later share
the same parameters: $\kappa =0.799432,$ $\lambda =0.940489,$ $\sigma =4;$
here $\theta =0.260229.$ Figure 2b illustrates a more pronounced switching
curve with $\kappa =3/13,$ $\lambda =9,$ $\sigma =4$ and $\theta =0.175$.%
\begin{figure}[!htbp]
\centering
 \begin{subfigure}[]{0.47\textwidth}
  \centering
    \includegraphics[width=\textwidth]{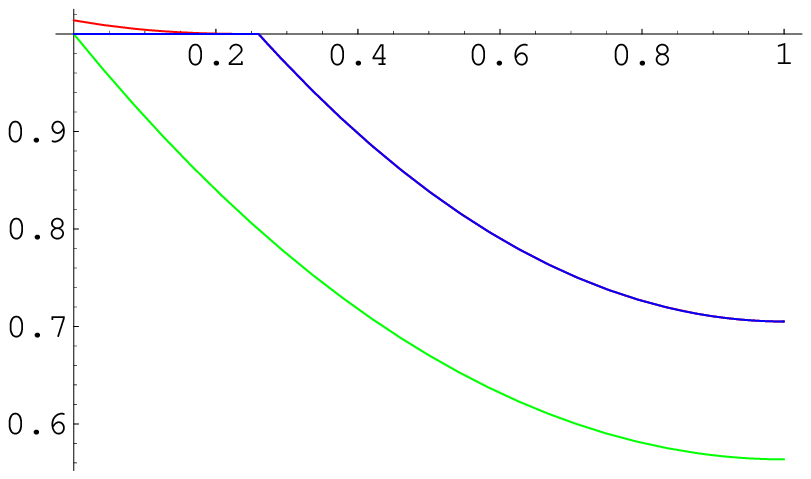}
    \caption*{Figure 2a. Switching curve $%
\protect\gamma _{t}^{\ast }$ red; sparing-first valuation $\protect\gamma %
_{t}$ blue, coalescing after $\protect\theta $; $\protect\gamma _{t}^{1}$
green.}
      \end{subfigure}
      \hfill
 \begin{subfigure}[]{0.47\textwidth}
  \centering
    \includegraphics[width=\textwidth]{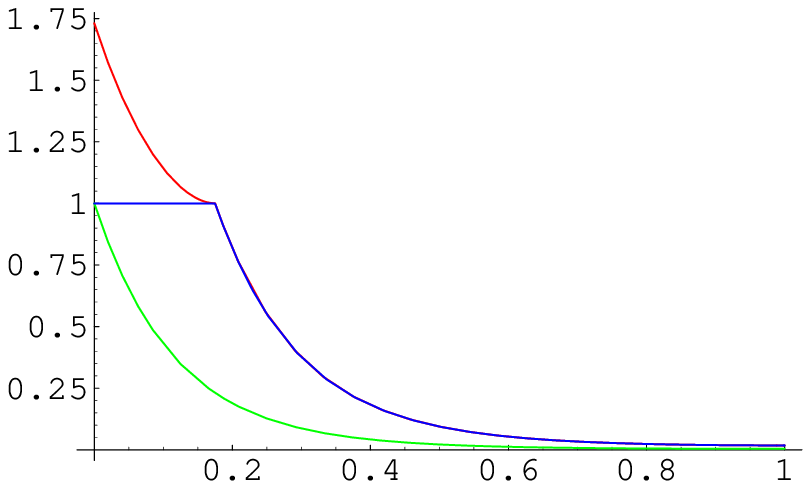}
    \caption*{Figure 2b. $\protect\gamma _{t}^{\ast }$ red; valuation $%
\protect\gamma _{t}$ blue; sparing-throughout valuation $\protect\gamma %
_{t}^{1}$ green; $\protect%
\gamma _{t}$ blue.}
      \end{subfigure}
\end{figure}

\textbf{Commentary to Example 1. }Here under silence, initially the
switching curve $\gamma _{t}^{\ast }$ (shown in red) is above the starting
firm value $\gamma _{0}=1$ and so candour ($\pi =0$) is initially optimal;
it is rational to infer that silence here means managers have received no
information, hence the valuation remains unchanged, until the switching time
is reached, as signified by the stationarity of the switching curve.
Thereafter, $\gamma _{t}^{\ast }$ falls below $\gamma =1$ and so the optimal
strategy yields a superior valuation to that given by $\gamma _{t}^{1}$
(which would have resulted from a policy of being sparing-throughout, i.e. $%
\pi \equiv 1$); here the valuation ignores the kind of silence that hides
bad news. In this time interval $\gamma _{t}^{\ast }$ and $\gamma _{t}$
coalesce, as predicted by the Non-mixing Theorem 1S. With a higher
intensity-value $\lambda $ of the private managerial news-arrival, the
switching time would come earlier, thus absorbing the higher chances of
ensuing privately received bad news, which strategically wants to be
withheld. The figures above graphically depict the reputational benefits to
a firm following a candid-first strategy. Starting at $t=0$ investors do not
downgrade the value of the firm when they see no disclosure, since they
infer this follows from non-observation of updated information. The blue
curve in both figures remains flat. This is in contrast to how investors
treat a firm applying a sparing strategy, for which they continually
downgrade firm value in the presence of silence. Thus the distance between
the blue and green lines reflects the reputational benefit of following a
candid strategy at times. In summary, the reputational benefit is in the
fact that investors do not downgrade firm value quite so heavily in the
presence of continuing silence.

\subsection{\textbf{Example 2}}
Here $\pi =1$ initially (initially sparing behaviour). In Figure 3, with $%
\kappa ,\lambda $ and $\sigma $ values as in Fig 2a, the switching time is $%
\theta =0.565997$.%
\begin{figure}[!htbp]
\begin{center}
\includegraphics[height=3.6cm]{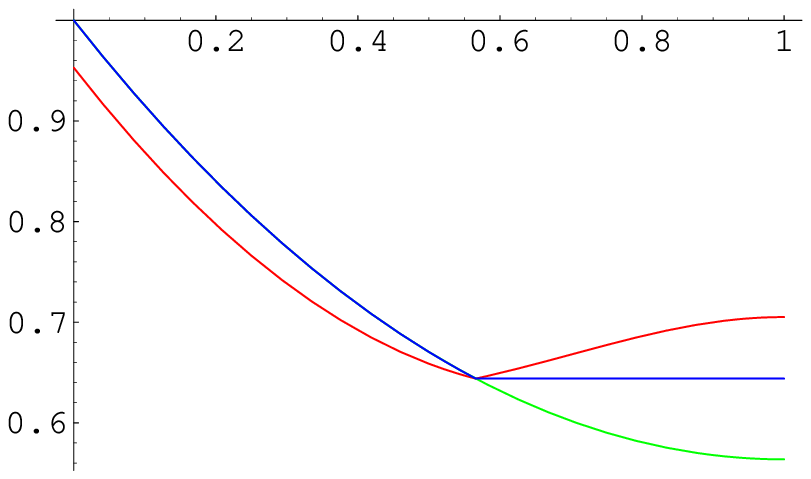}
\par
{Figure 3. $\protect%
\lambda$ and $\protect\kappa$ values as in Fig. 2a, but initially sparing.}
\end{center}
\end{figure}

\textbf{Commentary to Example 2. }Here, under silence, initially the
switching curve $\gamma _{t}^{\ast }$ (red) is below the firm valuation $%
\gamma _{t},$ albeit close, and $\gamma _{t}=\gamma _{t}^{1}$ is
consequently the dynamic disclosure-threshold (sparing policy threshold)
curve. Thereafter, $\gamma _{t}^{\ast }$ is above the $\gamma _{t}^{1}$
(green) curve, so it is optimal to switch to candour, which yields a
constant equilibrium valuation under silence (shown in blue). The valuation
subsequently omits to account for the kind of silence that hides bad news.
After the disclosure policy switch from $\pi =1$ to $\pi =0,$ it is rational
to infer that silence means managers have received no new information. With
a higher intensity $\lambda $ of private managerial news-arrival, the
switching time would come later, thus absorbing the higher chances of bad
news arrival needing strategically to be withheld.
\subsection{\textbf{Example 3 (non-constant }$\protect\kappa $\textbf{) }}
\begin{figure}[!htbp]
\begin{center}
\includegraphics[height=3.6cm]{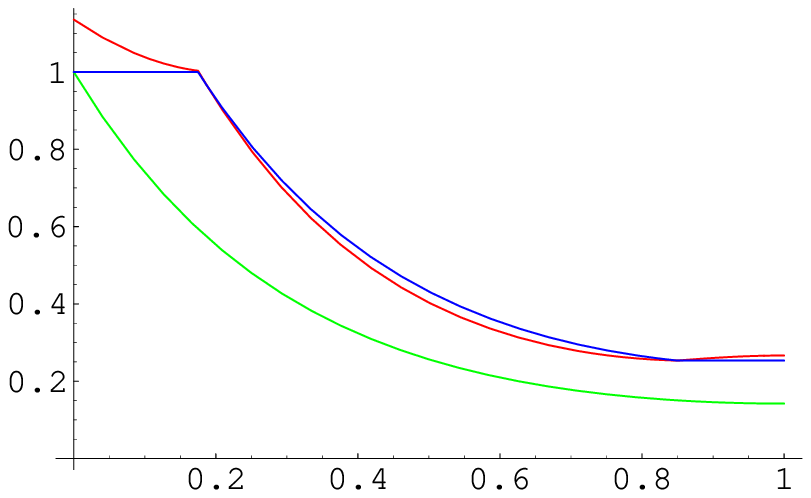}
\par
{Figure 4. Switching curve $%
\protect\gamma _{t}^{\ast }$ red; $\protect\gamma _{t}^{1}$ green; $\protect%
\gamma _{t}$ blue.}
\end{center}
\end{figure}

Here $\pi =0$ on $[0,\theta _{1})$ with $\theta _{1}=0.175$ and $\kappa
=0.594;$followed by $\pi =1$ on $[\theta _{1},\theta _{2})$ with $\theta
_{2}=0.85$ and $\kappa =0.533$ and finally by $\pi =0$ on $[\theta _{2},1].$
Throughout $\lambda =3.2$ and $\sigma =2.$

A choice of piecewise constancy may at first sight seem specious.
However, this is an arrangement of a pre-determined managerial reward
capable of being agreed by the shareholders, as the switching times are not
dynamically selected. Moreover, our analysis with constant $\kappa $ can be
adapted (by reference to the first mean-value theorem of
integration) to the general case of continuous $\kappa _{t}$ by replacing
within any inter-switching interval a proposed varying $\kappa _{t}$ by
some appropriate constant value (along the lines of a `certainty equivalent'
relative to the Poisson jumps), that value being intermediate between those
taken by $\kappa _{t}$ on that interval.

\section{Proofs}

The proof of Theorem 1 (actually in a stronger form) is in Appendix A.
Here we consider Theorems 2a and 2b. We recall that below $\alpha _{t}$ and $%
\beta _{t}$ are assumed constant. We begin with some preliminary
observations.

\bigskip

\noindent \textbf{Proposition 2.} (i) \textit{On any interval }$[\theta
^{\prime },\theta ]$\textit{\ where }$\pi =1,$\textit{\ the co-state
equation and solution take the form:}%
\begin{eqnarray*}
\frac{\mu _{t}^{\prime }}{\mu _{t}} &=&\lambda h(t)>0,\text{ for }\theta
^{\prime }<t<\theta , \\
\mu _{t} &=&K\exp [-\lambda (g(\theta )-g(t))],
\end{eqnarray*}%
\textit{so that }$\mu _{t},$\textit{\ being negative, is decreasing with }$%
K<0$\textit{\ and}%
\[
\partial \mu _{t}/\partial \lambda =(g(\theta )-g(t))(-\mu _{t})>0,\text{
for }\theta ^{\prime }<t<\theta .
\]

(ii) \textit{On any interval where }$\pi =0,$%
\[
\mu _{t}^{\prime }=(\beta _{t}-\alpha _{t})>0:\mu _{t}=-(\beta -\alpha )(K-t)%
\text{ if }\beta _{t}\equiv \beta \text{ and }\alpha _{t}\equiv \alpha \text{
with }K\leq 1.
\]%
\textit{Thus here }$\mu _{t}$\textit{\ is increasing.}

\textit{In particular, in both cases }$\mu _{t}$\textit{\ is either
non-constant or zero.}

\textit{\bigskip }

\noindent \textbf{Proof.} Since the co-state equation asserts that%
\[
\mu _{t}^{\prime }-\mu _{t}\pi _{t}\mathbb{\lambda }h(t)=(\beta _{t}-\alpha
_{t})(1-\pi _{t}),
\]%
the conclusions are immediate from the form of the differential equation.
\hfill $\square $

\bigskip

\noindent \textbf{Remark (Behaviour of }$\gamma $\textbf{). }If switches
occur at the three times $\theta _{1}<\theta _{2}<\theta _{3}$ with $\pi =1$
on $(\theta _{1},\theta _{2}),$ then $\gamma _{t}=e^{\lambda g(\theta
_{1})}e^{-\lambda g(t)}\gamma _{\theta _{1}}$ for $t\in \lbrack \theta
_{1},\theta _{2}],$ so
\[
\gamma _{\theta _{3}}=\gamma _{\theta _{2}}=\gamma _{\theta _{1}}e^{\lambda
g(\theta _{1})}e^{-\lambda g(\theta _{2})}.
\]%
Applying the formula inductively, if $\pi =0$ near $t=0,$ so that $\gamma
_{\theta _{1}}=1,$ then%
\[
\gamma _{\theta _{2n}}=e^{\lambda g(\theta _{1})}e^{-\lambda g(\theta
_{2})}e^{\lambda g(\theta _{3})}e^{-\lambda g(\theta _{4})}...e^{\lambda
g(\theta _{n-1})}e^{-\lambda g(\theta _{n})};
\]%
likewise if $\pi =1$ near $t=0,$ so that $\gamma _{\theta _{1}}=\gamma
_{\theta _{2}}=e^{-\lambda g(\theta _{1})},$ then%
\[
\gamma _{\theta _{2n}}=e^{-\lambda g(\theta _{1})}e^{\lambda g(\theta
_{2})}e^{-\lambda g(\theta _{3})}e^{\lambda g(\theta _{4})}...e^{\lambda
g(\theta _{n-1})}e^{-\lambda g(\theta _{n})}.
\]

\noindent \textbf{Corollary 1 (Final switching conditions).} \textit{If the
last two intervals of }$\pi $ \textit{constancy are given by }$\pi =0$%
\textit{\ switching at }$\theta $\textit{\ to }$\pi =1,$\textit{\ then }$\mu
_{t}=0$\textit{\ on }$[\theta ,1]$\textit{\ and near and to the left of }$%
\theta :$%
\[
\mu _{t}=(\beta -\alpha )(t-\theta ).
\]%
\textit{For the reversed strategy, if the last two aforementioned intervals
are given by }$\pi =1$\textit{\ switching at }$\theta $\textit{\ to }$\pi
=0, $\textit{\ then on }$[\theta ,1]$
\begin{eqnarray*}
\mu _{t} &=&(\beta -\alpha )(t-1),\text{ so that }\mu _{\theta }=(\beta
-\alpha )(\theta -1), \\
\mu _{t} &=&(\beta -\alpha )(\theta -1)e^{-\lambda \lbrack g(\theta )-g(t)]}%
\text{ for }t<\theta \text{ near }\theta .
\end{eqnarray*}%
\noindent \textbf{Corollary 2.} \textit{Being candid at all times (}$\pi
\equiv 0$\textit{) is optimal for }$\lambda $\textit{\ sufficiently low,
i.e. below a threshold depending on }$\kappa $ \textit{(equivalently,
depending on} $\alpha /\beta $\textit{).}

\bigskip

\noindent \textbf{Proof.} The assumption $\pi \equiv 0$ implies $\gamma
_{t}\equiv 1$ and $\gamma _{t}<\gamma _{t}^{\ast }.$ From Cor. 1, since $\mu
_{1}=0,$ we have $\mu _{t}=-(\beta -\alpha )(1-t),$ and so
\[
\gamma _{t}^{\ast }=\frac{e^{-\lambda g(t)}}{\kappa (1+\lambda h(t)(1-t))}>1%
\text{ iff }e^{-\lambda g(t)}>\kappa (1+\lambda h(t)(1-t)).
\]%
This holds for all $\lambda $ small enough (depending on $\kappa $), indeed,%
\[
\lim_{\lambda \rightarrow 0}\frac{e^{-\lambda g(1)}}{1+\lambda }=1>\kappa ,
\]%
so for all $\lambda $ small enough%
\[
e^{-\lambda g(1)}>\kappa (1+\lambda )
\]%
and, since $h(0)<1$ and $h(t)(1-t)$ decreases on $[0,1],$
\[
e^{-\lambda g(t)}>e^{-\lambda g(1)}>\kappa (1+\lambda )>\kappa (1+\lambda
h(t)(1-t)),
\]%
as $g(t)$ is increasing.\hfill $\square $

\bigskip

We need some details about the function $h$:

\bigskip

\noindent \textbf{Lemma 1.} \textit{The function} $h$ \textit{is decreasing
with }$h(1)=0$\textit{\ and }$h(t)<1$\textit{\ for} $t\in \lbrack 0,1].$

\bigskip

\noindent \textbf{Proof.} Recall that $h(t)=[2\Phi (\hat{\sigma}/2)-1],$
where $\hat{\sigma}=\sigma (1-t).$ So $h$ is decreasing (to $0)$, since%
\[
h^{\prime }(t)=-\sigma \exp \left( -[\sigma ^{2}(1-t)^{2}/8]\right) /\sqrt{%
2\pi }<0.
\]%
Note that $h(1)=0,$ as $\Phi (0)=1/2,$ and further that $h(0)>0$ as $\Phi
(\sigma /2)>1/2$ for $\sigma >0.$ Finally note that $h(t)<1,$ since $h(0)<1,$
the latter because%
\[
2\Phi (\sigma /2)-1<1\text{ as }\Phi (\sigma /2)<1.\qquad \qquad \qquad
\qquad \square
\]%
\hfill

If in equilibrium there is just one optimal switching $\theta $, it is
straightforward to derive its location property using a first-order
condition on $\theta $, as given in Theorems 2a and 2b above. The existence
conditions for the equilibria are derived from the general Hamiltonian
formulation; this requires the explicit derivation of the switching curve,
which is the content of a technical lemma and carries\ all the work for
theorems 2a (and likewise for 2b).

\bigskip

\noindent \textbf{Lemma 2a.}

\noindent \textit{If }$\pi _{t}=0$\textit{\ for }$t<\theta ,$\textit{\ and }$%
\pi _{t}=1$\textit{\ for }$t>\theta ,$\textit{\ is optimal, then }%
\[
\mu _{t}=\left\{
\begin{array}{cc}
\beta \kappa (t-\theta )=(\beta -\alpha )(t-\theta ), & t\leq \theta ,\text{
} \\
0, & t>\theta .%
\end{array}%
\right.
\]

\noindent \textit{Here the switching curve is given by}%
\[
\gamma _{t}^{\ast }=\left\{
\begin{array}{cc}
\left. \gamma _{t}^{1}\right/ \kappa \lbrack 1+\lambda (\theta -t)h(t)], &
t\leq \theta , \\
\left. \gamma _{t}^{1}\right/ \kappa , & t>\theta .%
\end{array}%
\right.
\]%
\textit{Here, }%
\[
\gamma _{t}^{\ast }=\gamma _{t},\text{ for }t>\theta .
\]

\noindent \textbf{Remark. }The graph of\textbf{\ }$\gamma _{t}^{\ast }$ is
stationary at $t=\theta ,$ for, writing $\partial _{t}$ for the time
derivative,%
\[
\kappa \partial _{t}\gamma _{t}^{\ast }=\left. \frac{-\lambda
h(t)e^{-\lambda g(t)}}{1+\lambda (\theta -t)h(t)}-\frac{-\lambda
h(t)+\lambda (\theta -t)h^{\prime }(t)}{[1+\lambda (\theta -t)h(t)]^{2}}%
e^{-\lambda g(t)}\right\vert _{t=\theta }=0.
\]

\noindent \textbf{Proof of Lemma 2a. }We will deduce $\mu _{t}$ directly
from the formal solution (Appendix A) via the integrating factor $\varphi $
there, which reduces to%
\[
\varphi (t)=\left\{
\begin{array}{cc}
1\text{ }, & t\leq \theta , \\
e^{-\lambda \lbrack g(t)-g(\theta )]}, & \theta <t<1.%
\end{array}%
\right.
\]%
For $t>\theta ,$ $\mu _{t}=0.$ For $t<\theta ,$ we have%
\[
\int_{t}^{\theta }\varphi _{s}\varphi _{t}^{-1}(\beta -\alpha )(1-\pi
_{s})_{=1}\text{ }\mathrm{d}s+\int_{\theta }^{1}\varphi _{s}\varphi
_{t}^{-1}(\beta -\alpha )(1-\pi _{s})_{=0}\text{ }\mathrm{d}s=-\mu _{t},
\]%
so%
\[
(\beta -\alpha )(\theta -t)=-\mu _{t}.
\]%
For consistency we need at $t=\theta =\theta (\lambda )$ that%
\[
\gamma _{\theta }^{\ast }=\frac{\gamma _{\theta }^{1}}{\kappa }=1:\kappa
=\gamma _{\theta }^{1}=e^{-\lambda g(\theta )}=\exp \left( -\lambda
\int_{0}^{\theta (\lambda )}h(s)\text{ }\mathrm{d}s\right) .
\]

As regards the coalescence, note that if $\theta $ solves $e^{-\lambda
g(\theta )}=\kappa $, as above, then, for $t\geq \theta ,$%
\[
\gamma _{t}=e^{-\lambda \lbrack g(t)-g(\theta )]}=\left. \gamma
_{t}^{1}\right/ \kappa =\left. e^{-\lambda g(t)}\right/ \kappa ,
\]%
rather than the expected inequality (\TEXTsymbol{>}); so as in Theorem 1S
(see Appendix A),
\[
\gamma _{t}^{\ast }=\gamma _{t},
\]%
and since $\gamma _{t}^{\prime }=-\lambda h(t)\gamma _{t},$ it follows that $%
\pi =1$.\hfill $\square $

\bigskip

Armed with Lemma 2a we proceed to

\bigskip

\noindent \textbf{Proof of Theorem 2a.} We begin with the location
condition. Since $\pi _{t}=0$ for $t\in $ $[0,\theta ),$ then $\gamma
_{t}\equiv 1$ and $\gamma _{t}^{1}=e^{-g(t)}$ on this interval and so the
objective function reduces to%
\[
\int_{0}^{\theta }e^{-\lambda g(t)}dt-\kappa \theta ,
\]%
since $1-\pi _{t}=0$ on $[\theta ,1]$, which yields a zero contribution.
Differentiation w.r.t. $\theta $ yields the first-order condition%
\[
e^{-\lambda g(\theta )}-\kappa =0,
\]%
which on re-arrangement yields the claim.

We turn to the existence condition. To ensure that $\gamma _{0}^{\ast
}>\gamma _{0}$ requires by Lemma 2a that
\[
\left. 1\right/ \kappa \lbrack 1+\lambda \theta h(0)]>1,\text{ equivalently }%
\lambda <(\kappa ^{-1}-1)\left/ (\theta h(0))\right. .
\]%
(The upper bound on $\lambda $ is illustrated by the green curve in Figure 5
below.) This also guarantees that $\gamma _{t}^{\ast }>\gamma _{t}$ for $%
t<\theta ,$ since $(\theta -t)h(t)$ is decreasing. It is also required that $%
\gamma _{\theta }^{\ast }=1,$ i.e. $e^{-\lambda g(\theta )}=\kappa ,$ and so
this holds iff%
\[
-\log \kappa /g(\theta )=\lambda <(\kappa ^{-1}-1)\left/ (\theta
h(0))\right. .
\]%
(The form on the left side here is illustrated by the red curve in Figure
5.) By Lemma 2a, $\gamma _{t}^{\ast }=\gamma _{t}$ for $t>\theta .$

The displayed inequality is feasible iff%
\[
\log (\kappa ^{-1})/(\kappa ^{-1}-1)<g(\theta )\left/ (\theta h(0))\right. .
\]

Finally we compute the rate of change of $\theta $ w.r.t. $\lambda $ from
the location condition expressed as%
\begin{eqnarray*}
\lambda \int_{0}^{\theta (\lambda )}h(s)\text{ }\mathrm{d}s &=&-\log \kappa
>0: \\
\int_{0}^{\theta (\lambda )}h(s)\text{ }\mathrm{d}s+\lambda h(\theta
(\lambda ))\theta ^{\prime }(\lambda ) &=&0,\text{ so }\theta ^{\prime
}(\lambda )=\frac{\log \kappa }{\lambda ^{2}h(\theta (\lambda ))}<0.
\end{eqnarray*}%
So the larger is $\lambda $ the smaller is $\theta .$\hfill $\square $

\bigskip

From where the horizontal blue line $\lambda =10$ in the figure intersects
the green curve one may drop vertically to the red curve to obtain a value
of $\lambda $ which lies below the green and on the red curve and lower blue
line.%

\begin{figure}[tbp]
\begin{center}
\includegraphics{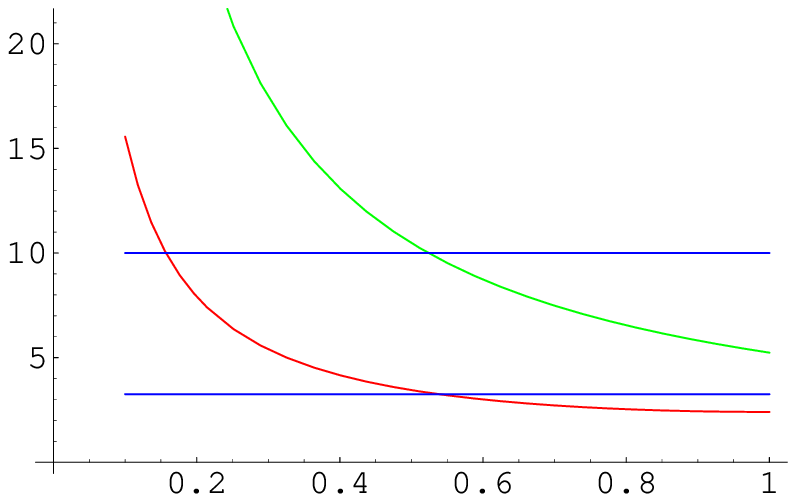}
\par
{Figure 5. $\lambda \,$%
 delimited above by
$(\kappa ^{-1}-1)\left/ (\theta h(0))\right.$%
and below by $\log \kappa^{-1} /g(\theta )$ against $\theta .$
}
\end{center}
\end{figure}

\bigskip

We again begin with a technical Lemma.

\bigskip

\noindent \textbf{Lemma 2b. }\textit{Assume that }$\alpha _{t},\beta _{t}$
\textit{are constant.}

\noindent \textit{In an equilibrium where }$\pi _{t}=1$\textit{\ for }$t<$%
\textit{\ }$\theta $\textit{\ and }$\pi _{t}=0$\textit{\ for }$t>\theta ,$%
\textit{\ one has}%
\[
\mu _{t}=\left\{
\begin{array}{cc}
(\beta -\alpha )(\theta -1)[\exp -\lambda \int_{t}^{\theta }h(u)\text{
\textrm{d}}u], & t\leq \theta , \\
(\beta -\alpha )(t-1), & \theta <t<1.%
\end{array}%
\right.
\]

\noindent \textit{Here}%
\[
\gamma _{t}^{\ast }=\left\{
\begin{array}{cc}
\left. \gamma _{t}^{1}\right/ \kappa \lbrack 1+\lambda (1-\theta )h(t)\exp
[-\lambda (g(\theta )-g(t))], & t<\theta ,\text{ } \\
\left. \gamma _{t}^{1}\right/ \kappa \lbrack 1+\lambda (1-t)h(t)], & t\geq
\theta .%
\end{array}%
\right.
\]%
\textit{This curve defines the switching time }$\theta $\textit{\ as the
intersection time }$t=\theta $\textit{\ of }$\gamma _{t}^{\ast }$ \textit{%
with }$\gamma _{t}^{1}$\textit{\ when}%
\[
\lambda (1-\theta )h(\theta )=\kappa ^{-1}-1,\text{i.e. }1=\theta +\frac{%
\kappa ^{-1}-1}{\lambda h(\theta )}.
\]%
\bigskip

\noindent \textbf{Proof.} As in Lemma 2a, we again deduce $\mu _{t}$
directly from the integrating factor (see Appendix A), which here is%
\[
\varphi (t)=\left\{
\begin{array}{cc}
e^{-\lambda g(t)}\text{ }, & t\leq \theta , \\
e^{-\lambda g(\theta )}, & \theta <t<1.%
\end{array}%
\right.
\]%
For $t>\theta ,$ as $s>\theta $ and $\pi =0$ below,%
\[
\int_{t}^{1}\varphi _{s}\varphi _{t}^{-1}(\beta _{s}-\alpha _{s})(1-\pi _{s})%
\text{ }\mathrm{d}s=-\mu _{t},\text{ i.e. }(\beta -\alpha )(1-t)=-\mu _{t}.
\]%
For $t<\theta $ we have%
\begin{eqnarray*}
-\mu _{t} &=&\int_{t}^{\theta }\varphi _{s}\varphi _{t}^{-1}(\beta
_{s}-\alpha _{s})(1-\pi _{s})_{=0}\text{ }\mathrm{d}s+\int_{\theta
}^{1}\varphi _{s}\varphi _{t}^{-1}(\beta _{s}-\alpha _{s})(1-\pi _{s})_{=1}%
\text{ }\mathrm{d}s \\
&=&[\exp \lambda \int_{0}^{t}h(u)\text{ }\mathrm{d}su][\exp -\lambda
\int_{0}^{\theta }h(u)\text{ }\mathrm{d}u](\beta -\alpha )(1-\theta ), \\
&=&(\beta -\alpha )(1-\theta )[\exp -\lambda \int_{t}^{\theta }h(u)\text{ }%
\mathrm{d}u].\qquad \qquad \qquad \qquad \square
\end{eqnarray*}

We may now prove Theorem 2b.

\bigskip

\noindent \textbf{Proof of Theorem 2b. }We begin with the location
condition. Since $\pi _{t}=1$ on $[0,\theta )$, $\gamma _{t}=\gamma
_{t}^{1}\equiv e^{-g(t)}$ on $[0,\theta )$. Here the objective function
reduces to%
\[
\int_{\theta }^{1}e^{-\lambda g(t)}dt-\kappa e^{-\lambda g(\theta
)}(1-\theta ),
\]%
since there is a zero contribution to the objective function on $[0,\theta
]. $ Differentiation w.r.t. $\theta $ yields the first-order condition%
\[
-e^{-\lambda g(\theta )}+\kappa e^{-\lambda g(\theta )}+\kappa e^{-\lambda
g(\theta )}\lambda h(\theta )(1-\theta )=0,
\]%
which on re-arrangement yields the claim.

We turn to the existence condition. By Lemma 2b, intersection at $t=\theta $
of $\gamma _{t}^{\ast }$ with $\gamma _{t}^{1}$ occurs iff%
\[
\lambda (1-\theta )h(\theta )\exp [-\lambda \int_{\theta }^{\theta }h(u)%
\text{ \textrm{d}}u]=\lambda (1-\theta )h(\theta )=\kappa ^{-1}-1.
\]%
Combining this with the requirement that $\gamma _{0}^{\ast }<1$ yields%
\begin{eqnarray*}
\left. 1\right/ \kappa \lbrack 1+\lambda (1-\theta )h(0)\exp [-\lambda
g(\theta )] &<&1, \\
\lambda (1-\theta )h(\theta ) &=&\kappa ^{-1}-1<\lambda (1-\theta )h(0)\exp
[-\lambda g(\theta )].
\end{eqnarray*}%
This holds for some $\kappa $ iff%
\begin{equation}
h(\theta )/h(0)<\exp [-\lambda g(\theta )],\text{ i.e. }\lambda <-[\log
h(\theta )/h(0)]/g(\theta ),
\end{equation}%
yielding a bound on $\lambda $ in terms of the switching time $\theta $ (as
illustrated by a green graph in Figure 6 below). From here we obtain%
\begin{equation}
\kappa ^{-1}-1=\lambda (1-\theta )h(\theta )<(1-\theta )h(\theta )(-[\log
h(\theta )/h(0)]/g(\theta )),
\end{equation}%
in turn a lower bound on $\theta $ (as illustrated by a red graph in Figure
6).

Finally we compute the rate of change of $\theta $ w.r.t. $\lambda $ from
the location condition. Here for $\theta =\theta (\lambda )$ we have, as $%
h^{\prime }(\theta )<0,$%
\[
\lambda (1-\theta )h(\theta )=\frac{\alpha }{\beta -\alpha },\text{ so }%
\theta ^{\prime }(\lambda )=\frac{(1-\theta )h(\theta )}{\lambda \lbrack
h(\theta )-(1-\theta )h^{\prime }(\theta )]}>0.\qquad \qquad \square
\]%

\begin{figure}[!htbp]
\begin{center}
\includegraphics{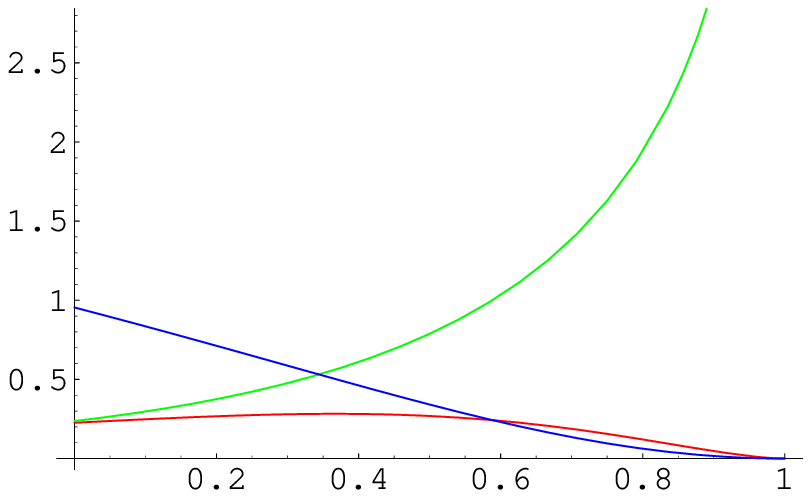}
\par
{Figure 6. Bounds for selecting $%
\protect\lambda .$ Graph in blue of $h(\protect\theta )$ against $\protect%
\theta .$}
\end{center}
\end{figure}

\noindent \textbf{Remark.} In Figure 6 above the red curve traces possible
values of the function in the first condition (1) above (and in Theorem 2b\
earlier); the green curve corresponds to the function in the second display
(2) above. The blue curve identifies the $\theta $ value given a horizontal $%
\lambda $ value. Thus $\lambda $ must lie on the portion of the blue curve
lying in between the red and green.

\section{Double and multiple-switching}

This section continues to study multi-switching, as a complement to our main theme
that single-switching should be regarded as the preferred equilibrium type. We saw this is the case only for a {\textit{candid-first}} equilibrium. Below we explain the numerical evidence for the absence of a second switch in a {\textit{sparing-first}} equilibrium.

We begin with a useful observation which refers to an equilibrium condition that must hold at an optimized switching location.
\bigskip

\noindent \textbf{Lemma 3. }\textit{In an equilibrium with two optimal
switching points }$0<\theta _{1}<\theta _{2}<1$\textit{\ for which }$\gamma
_{t}^{\ast }<\gamma _{t}$\textit{\ on }$[0,\theta _{1}),$\textit{\ where }$%
\pi =1,$ \textit{it is the case that}%
\[
\lambda <\eta (\sigma ):=-h^{\prime }(0)/h(0)^{2}.
\]

\noindent \textbf{Proof.} See Appendix B. \hfill $\square $.\bigskip

\noindent \textbf{Initial sparing policy: numerical evidence for the absence
of a second switching point.} The existence of a second optimal switching
time for a sparing-first equilibrium requires by Theorem 3b that $\theta
=\theta _{1},$ for $\theta _{1}$ the first switching point, solve%
\[
T(\theta ):=g(\theta +\frac{1-\kappa }{\lambda \kappa h(\theta )})-g(\theta
)+\frac{1}{\lambda }\log \kappa =0,
\]%
since $g(\theta_0)=g(0)=0$; this would yield a second switching point at $\theta =\theta _{2}$ given by%
\[
\theta _{2}=\theta _{1}+\frac{1-\kappa }{\lambda \kappa h(\theta_{1})}.
\]%
By Proposition A below, this requires either $\lambda $ to be small enough or $%
\kappa $ large enough (otherwise  $T(\theta)>0$). By Lemma 3, we may study the solvability of the
equation $T(\theta )=0$ by taking $\lambda =\lambda _{\text{max}}=\eta
(\sigma )=-h^{\prime }(0)/h(0)^{2}$ and graphing against the aggregate volatility $\sigma $ the
expression%
\[
T(0)=g\left( \frac{\kappa ^{-1}-1}{h(0)\lambda _{\text{max}}}\right) +\frac{1%
}{\lambda _{\text{max}}}\log \kappa ,
\]%
with $\kappa $ near unity, say at $0.99$, to determine whether or not $%
T(0)>0.$ It turns out that $T(0)<0$ on a wide range of $\sigma .$
Furthermore, the sign of%
\[
T^{\prime }(0)=h\left( \frac{\kappa ^{-1}-1}{h(0)\lambda _{\text{max}}}%
\right) \left( 1+\frac{-h^{\prime }(0)}{h(0)^{2}}\frac{\kappa ^{-1}-1}{%
\lambda _{\text{max}}}\right) -h(0)
\]%
is also found to be negative on a wide range of $\sigma .$ Thus numerical
evidence supports the hypothesis that $T(\theta )$ remains negative and
bounded away from $0$, so that the equation $T(\theta )=0$ is insolvable, and
hence there does not exist a second switching point. This implies that for $%
\kappa _{t}$ constant there is at most one switching point in an equilibrium
which initially implements a sparing policy.

\bigskip

{\textbf{Remark.}} In general, how many switches should there be in equilibrium? This question
can be resolved by reference to the preferences either of the managers or
alternatively of a representative investor in the market. The former chooses
the policy with higher reward to the managers; here, the pay-for-performance
sides with multiple switching, since it boosts the valuation and so also the
managerial payoff. Resolution by reference to investors requires the
introduction of a notional cost of complexity. But the underlying assumption
must remain that the market has found its way to an appropriate equilibrium -- cf. \S 2.

\bigskip

We return to examining the formula
\[
T_{\kappa \lambda }(\theta ):=\theta +\frac{\kappa ^{-1}-1}{\lambda h(\theta
)},
\]
which defines for any optimal switching point $\theta $ the linked (coupled)
successor switching point in terms of model parameters ($\kappa ,\lambda $
and $\sigma ,$ the latter embedded in $h$). The constraint $T_{\kappa
\lambda }(\theta )\leq 1$ holds iff%
\begin{eqnarray*}
\kappa ^{-1} &\leq &(1+\lambda (1-\theta )h(\theta ))\leq 1+\lambda h(0), \\
\kappa  &\geq &(1+\lambda h(0))^{-1},
\end{eqnarray*}%
placing a lower bound on $\kappa $ in terms of $\lambda $ (and $\sigma $).
This implies that, for given $\sigma $ and $%
\lambda ,$ there is a critical value $\bar{\kappa},$ such that the key
equation above is insolvable for $\kappa <\bar{\kappa}.$ It is also true
that, for all small enough $\kappa ,$ there is a critical value of $\lambda ,
$ above which the key equation $T(\theta )=0$ is not solvable, as clarified
by the next result.

\bigskip

\noindent \textbf{Proposition A.} \textit{The following inequality holds}%
\[
g(\theta )+\frac{\log \kappa ^{-1}}{\lambda }<g\left( \theta +\frac{\kappa
^{-1}-1}{\lambda h(\theta )}\right) \text{ for }T_{\kappa, \lambda }(\theta )\leq 1,%
\text{ }
\]%
\textit{provided }$\kappa $\textit{\ and }$\lambda $\textit{\ satisfy }%
\[
\frac{\log \kappa ^{-1}}{\kappa ^{-1}-1}\leq \frac{1}{2}\text{ and }\lambda
>(\kappa ^{-1}-1)\cdot \mathrm{\max }\left\{ \frac{|h^{\prime }(1)|}{h(\bar{%
\theta})^{2}},\frac{1}{h(0)}\right\} .
\]%
\noindent \textbf{Proof. }See Appendix B. \hfill $\square $.

\bigskip

\noindent \textbf{Remark.} To see why this last result holds, note that
expansion of $g$ near $\theta _{1}$ gives a first approximation to the
right-hand side of%
\[
g(\theta _{1})+\frac{\kappa ^{-1}-1}{\lambda h(\theta _{1})}g^{\prime
}(\theta _{1})=g(\theta _{1})+\frac{1-\kappa }{\lambda \kappa }.
\]%
Using the approximation $\log (1+x)\simeq x$ for $x$ small with $x=\kappa -1$
yields%
\[
-\frac{1}{\lambda }\log \kappa \simeq \frac{1-\kappa }{\lambda }.
\]%
This shows both sides can be close for $\kappa $ close to $1,$ in view of
the comparison%
\[
\frac{1-\kappa }{\lambda }\text{ versus }\frac{1-\kappa }{\lambda \kappa }.
\]%
However, for small $\kappa ,$ $\frac{1-\kappa }{\lambda \kappa }$ is too far
from $\frac{1-\kappa }{\lambda }$ and so even double-switching is then not
possible.

\bigskip

We close this section by noting a complementary equilibrating result
linking the switching point of Theorem 3a to the prior switching point.
Corollary 1 in \S 5 implies that the switching curve is made to fall on
intervals where $\pi _{t}=0$ by a corresponding linear decreasing effect of $%
|\mu _{t}|,$ helped by the decreasing nature of the instantaneous
protective-put, $h(t)$. Contrarily, it is made to rise on intervals where $%
\pi _{t}=1$ by the ensuing exponential nature, $e^{\lambda g(t)},$ embedded
in $|\mu _{t}|$. It is not therefore surprising to see $h(t)e^{\lambda g(t)}$
as the \textit{equilibrating} factor in (equil) below. Since $%
h(t)e^{\lambda g(t)}$ has a maximum at some $t\in (\theta _{i-1},\theta
_{i})$ satisfying%
\[
\lambda =-h^{\prime }(t)/h(t)^{2}
\]%
(cf. Lemma 3), Proposition B below determines $\theta _{i}$ from $\theta
_{i-1}$ when $\pi =1$ on $[\theta _{i-1},\theta _{i}).$

\bigskip

\noindent \textbf{Proposition B.} \textit{As in the setting of Theorem 3a, for each switching point }$\theta _{i}$\textit{\ that is a right endpoint
of an interval with }$\pi =1$ \textit{constancy}, \textit{the preceding
switching point satisfies the equilibrating relation}%
\begin{equation}
h(\theta _{i-1})e^{\lambda g(\theta _{i-1})}=h(\theta _{i})e^{\lambda
g(\theta _{i})}.  \tag{equil}
\end{equation}

\bigskip

\noindent \textbf{Proof. }See Appendix B.\hfill $\square $

\bigskip

\section{Conclusions}

The disclosure model of Dye [Dye] alerts us to consider the implications of
firms remaining silent between mandatory disclosure dates. The model
developed here shows why management of a firm may benefit from establishing
a reputation for being candid on some time intervals, voluntarily disclosing
all news, good or bad. At issue is when would one expect to see such
behaviour in an equilibrium and whether it is likely to be time-invariant,
once established. Corollary 2 (Section\ 5) establishes that if the news
intensity-arrival rate is sufficiently low an equilibrium exists in which
managers are always candid. The comparative statics\ of $\theta (\lambda )$
in Theorem 2a establishes that, as the news-arrival rate rises, eventually
the optimal policy for management is to switch to a sparing disclosure
policy; a similar effect is caused by the remaining factors in the model,
namely of time-to-expiry (to the next mandatory disclosure time) and
pay-for-performance ratio $\kappa .$ In the model, reputation for adoption
of a candid disclosure strategy is derived endogenously and we see that for
higher levels of news-arrival such a strategy will not be time-invariant --
managers will start to `burn their reputation' (switching from candid
reporting to sparing disclosure) the closer they get to a mandatory
disclosure date. If the aim is to understand asset pricing in a
continuous-time setting, then this model provides insights into how firm
management will voluntarily disclose information to update markets in
between mandatory disclosure dates. Litigation concerns may ensure that very
negative news is always disclosed; nevertheless, as this model shows, once
management switch out of candid disclosure into a sparing policy they will
tend to hide "slightly" negative news, both when close to a mandatory
disclosure time and when their private news-arrival rate is higher.

\section{References}

\noindent \lbrack AchMK] V. Acharya, P. DeMarzo, and I. Kremer.: Endogenous
information flows and the clustering of announcements. Amer. Econ. Rev.
101.7, 2955-79 (2011)

\noindent \lbrack BerMTV] J. Bertomeu, I. Marinovic, S. J. Terry, and F.
Varas.: The Dynamics of Concealment, J. Fin. Econ. 143.1 (2022), 227-246.

\noindent \lbrack BeyD] A. Beyer and R. A. Dye.: Reputation management and
the disclosure of earnings forecasts, Review of Accounting Studies 17,
877--912 (2012)

\noindent \lbrack BinK] N. H. Bingham and R. Kiesel.: Risk-neutral
valuation. Springer (1998)

\noindent \lbrack BreP] A. Bressan and B. Piccoli.: Introduction to
Mathematical Control Theory, AIMS\ (2007)

\noindent \lbrack Dav84] M. H. A. Davis.: Piecewise-deterministic Markov
processes: A general class of non-diffusion stochastic models, J. R. Stat.
Soc., Ser. B, 46.3, 353-388 (1984)

\noindent \lbrack Dav93] M. H. A. Davis.: Markov models and optimization.
Monographs on Statistics and Applied Probability 49, Chapman \& Hall, (1993)

\noindent \lbrack Dye] R. A. Dye.: Disclosure of Nonproprietary Information.
Journal of Accounting Research 23, 123-145 (1985)

\noindent \lbrack EinZ] E. Einhorn and A. Ziv.: Intertemporal Dynamics of
Corporate Voluntary Disclosures. J. Acc. Res. 46.3, 567-589 (2008)

\noindent \lbrack FudT] D. Fudenberg and J. Tirole.: Game Theory, MIT (1991)

\noindent \lbrack GieO] M. B. Gietzmann and A.J. Ostaszewski.: The sound of
silence: equilibrium filtering and optimal censoring in financial markets.
Adv. in Appl. Probab. 48A, 119--144 (2016)

\noindent \lbrack GieOS] M. B. Gietzmann, A. J. Ostaszewski and M. H. G. Schr%
\"{o}der.: Guiding the guiders: foundations of a market-driven theory of
disclosure. Stochastic modeling and control (2020), 107--132, Banach Center
Publ., 122pp.

\noindent \lbrack Gru] M. D. Grubb.: Developing a Reputation for Reticence.
J. Econ. Man. Strat. 20.1, 225-268 (2011)

\noindent \lbrack GutKS] I. Guttman, I. Kremer and A. Skrzypacz. "Not only
what but also when: A theory of dynamic voluntary disclosure." Amer. Econ.
Rev. 104.8, 2400-2420 (2014)

\noindent \lbrack KamG] E. Kamenica and M. Gentzkow.: Bayesian Persuasion.
Amer. Econ. Rev. 1010.6, 2590-2615 (2011)

\noindent \lbrack Lib] D. Liberzone.: Calculus of Variations and Optimal
Control Theory: A Concise Introduction, Princeton (2012)

\noindent \lbrack Lok] A. L\o kka.: Detection of disorder before an
observable event, Stochastics 79, 219--231 (2007)

\noindent \lbrack MarV] I. Marinovic and F. Vargas.: No news is good news:
Voluntary disclosure in the face of litigation. The RAND Journal of
Economics 47.4, 822-856 (2016)

\noindent \lbrack OstG] A. J. Ostaszewski and M. B. Gietzmann.: Value
creation with Dye's disclosure option: optimal risk-shielding with an upper
tailed disclosure strategy. Review of Quantitative Finance and Accounting
31, 1--27 (2008)

\noindent \lbrack Sas] A. Sasane.: Optimization in Function Spaces, Dover
(2016)

\noindent \lbrack Tro] J. L. Troutman.: Variational Calculus and Optimal
Control, Springer (1996)

\noindent \lbrack Wei] J. W. Weibull.: Evolutionary Game Theory, MIT (1997)

\section*{Appendix A}

\noindent \textbf{Proposition 3.} $\mu _{t}\leq 0$\textit{\ for }$t\in
\lbrack 0,1].$

\bigskip

\noindent \textbf{Proof.} This follows from a formal solution the co-state
equation. With $\pi _{t}$ piecewise-constant and $h$ as in \S 3.3, take
\[
\varphi (t)=\exp \left( -\int_{0}^{t}\pi _{u}\mathbb{\lambda }h(u)\text{
\textrm{d}}u\right) \geq 0
\]%
(the integrating factor for the co-state equation), a decreasing function so
that $\varphi (t)\leq \varphi (s)$ for $s\geq t,$ as indeed,%
\[
\varphi (s)/\varphi (t)=\exp \left( -\int_{t}^{s}\pi _{u}\mathbb{\lambda }%
h(u)\text{ \textrm{d}}u\right) \leq 1.
\]%
With $\mu _{t}$ piece-wise smooth and $\mu _{1}=0$ (see \S 3.4), integration
from $t$ to $1$ yields%
\begin{eqnarray*}
\frac{d}{dt}\mu _{t}\varphi _{t} &=&\varphi _{t}(\beta _{t}-\alpha
_{t})(1-\pi _{t})\geq 0\text{ as }\beta _{t}\geq \alpha _{t}\text{ (see \S %
3.3)} \\
0-\mu _{t}\varphi _{t} &=&\int_{t}^{1}\varphi _{s}(\beta _{s}-\alpha
_{s})(1-\pi _{s})\text{ }\mathrm{d}s, \\
-\mu _{t} &=&\varphi _{t}^{-1}\int_{t}^{1}\varphi _{s}(\beta _{s}-\alpha
_{s})(1-\pi _{s})\text{ }\mathrm{d}s\geq 0, \\
-\mu _{t} &=&\int_{t}^{1}\varphi _{s}\varphi _{t}^{-1}(\beta _{s}-\alpha
_{s})(1-\pi _{s})\text{ }\mathrm{d}s\leq \int_{t}^{1}(\beta _{s}-\alpha
_{s})(1-\pi _{s})\text{ }\mathrm{d}s.
\end{eqnarray*}%
Above we used the blanket assumption that $\alpha _{t}\leq \beta _{t},$ so
we conclude that $\mu _{t}\leq 0.$\hfill $\square $

\bigskip

\noindent \textbf{Corollary 4. }$|\mu _{t}|$ \textit{is bounded on }$[0,1]$%
\textit{\ when }$\alpha _{t},\beta _{t}$ \textit{\ are constant. }

\bigskip

\noindent \textbf{Proof. }This follows again from the blanket assumption and
from
\begin{eqnarray*}
0 &\leq &-\mu _{t}=\int_{t}^{1}\varphi _{s}\varphi _{t}^{-1}(\beta
_{s}-\alpha _{s})(1-\pi _{s})\text{ }\mathrm{d}s\leq \int_{t}^{1}(\beta
_{s}-\alpha _{s})(1-\pi _{s})\text{ }\mathrm{d}s \\
&\leq &(\beta -\alpha )\int_{t}^{1}(1-\pi _{s})\text{ }\mathrm{d}s\text{ if }%
\alpha ,\beta \text{ are constant}\leq (\beta -\alpha )(1-t).\qquad \qquad
\square
\end{eqnarray*}%
\hfill

\bigskip

\noindent \textbf{Theorem 1S (Non-mixing Theorem -- Strong Form). }\textit{%
Assume }$\alpha _{t},\beta _{t}$ \textit{are constant.}\newline
(i) \textit{If the state trajectory and switching curve coalesce on an
interval }$I$\textit{, then }$\pi \equiv 1$ \textit{on }$I.$\newline
(ii) \textit{A mixing control with }$\pi _{t}\in (0,1)$\textit{\ is
non-optimal over any interval of time.}

\noindent \textbf{Proof. }If a mixing control occurs, then, from the
Hamiltonian maximisation, it follows that $\gamma _{t}=\gamma _{t}^{\ast }$
on an interval of time; hence (ii) follows from (i), by contradiction. To
prove (i), we compute in Step 1 the corresponding control $\pi _{t}$ from
the equilibrium equation and then, in Step 2, show that this control does
not satisfy the co-state equation unless $\mu _{t}\equiv 0,$ in which case $%
\pi _{t}\equiv 1.$

\textbf{Step 1.} Since $\gamma _{t}$ and so also $\gamma _{t}^{\ast }$
satisfies the equation (cont-eq-$\pi )$ on $I$,
\[
(\gamma _{t}^{\ast })^{\prime }=-\pi _{t}\lambda \gamma _{t}^{\ast }h(t).
\]%
For ease of calculations, write
\[
\gamma _{t}^{\ast }=\psi \gamma _{t}^{1}\text{ with }\psi =\psi _{t}:=\frac{1%
}{[\kappa -\mu _{t}\lambda h(t)/\beta ]}>0.
\]%
Then, as $(\gamma _{t}^{1})^{\prime }=-\lambda \gamma _{t}^{1}h,$%
\[
(\gamma _{t}^{\ast })^{\prime }=\psi ^{\prime }\gamma _{t}^{1}+\psi \gamma
_{t}^{\prime 1}=\psi ^{\prime }\gamma _{t}^{1}+\psi \lbrack -\lambda \gamma
_{t}^{1}h].
\]%
So, from the equation (cont-eq-$\pi ),$%
\begin{eqnarray*}
-\pi _{t}\lambda h &=&\frac{\gamma _{t}^{\prime \ast }}{\gamma _{t}^{\ast }}=%
\frac{\psi ^{\prime }\gamma _{t}^{1}-\psi \lambda h\gamma _{t}^{1}}{\psi
\gamma _{t}^{1}}=\frac{\psi ^{\prime }-\psi \lambda h}{\psi }, \\
\pi _{t} &=&\frac{\psi \lambda h-\psi ^{\prime }}{\psi \lambda h}=1-\frac{%
\psi ^{\prime }}{\psi \lambda h}.
\end{eqnarray*}%
Substituting for $\psi $ and $\psi ^{\prime }$ and writing $h_{t}$ for $h(t)$
yields
\begin{eqnarray*}
\pi _{t} &=&1-\frac{[\kappa -\mu _{t}\lambda h_{t}/\beta ]}{\lambda h_{t}}%
\frac{(\mu _{t}h_{t})^{\prime }\lambda /\beta }{[\kappa -\mu _{t}\lambda
h_{t}/\beta ]^{2}}=1-\frac{(\mu _{t}h_{t})^{\prime }/\beta }{h_{t}[\kappa
-\mu _{t}\lambda h_{t}/\beta ]} \\
&=&1-\frac{\psi }{\beta }\frac{(\mu _{t}h_{t})^{\prime }}{h_{t}}=1-\frac{%
\psi }{\beta }\left( \frac{\mu _{t}^{\prime }h_{t}}{h_{t}}+\frac{\mu
_{t}h_{t}^{\prime }}{h_{t}}\right) =1-\frac{\psi }{\beta }\left( \mu
^{\prime }+\mu \frac{h^{\prime }(t)}{h(t)}\right) .
\end{eqnarray*}

\textbf{Step 2. }We now substitute this value for $\pi _{t}$ into the
co-state equation,%
\[
-\mu _{t}^{\prime }+\mu _{t}(\mathbb{\lambda }\pi _{t}h(t))=-\beta \kappa
(1-\pi _{t}),
\]%
(as $\kappa =1-\alpha /\beta ).$ We compute $\pi _{t}$ from the co-state
equation to be%
\[
\pi _{t}=\frac{\mu _{t}^{\prime }-\beta \kappa }{[\mu _{t}\mathbb{\lambda }%
h(t)-\beta \kappa ]},
\]%
the division being valid, since $\mu _{t}\mathbb{\lambda }h(t)\leq 0$ and $%
\beta \kappa >0.$ So,
\[
\pi _{t}=\frac{\mu _{t}^{\prime }/\beta -\kappa }{[\mu _{t}\mathbb{\lambda }%
h(t)/\beta -\kappa ]}=-\psi _{t}[\mu _{t}^{\prime }/\beta -\kappa ]=\psi
_{t}[\kappa -\mu _{t}^{\prime }/\beta ].
\]%
Consistency of the two formulas requires that%
\[
\psi _{t}[\kappa -\mu _{t}^{\prime }/\beta ]=1-\frac{\psi _{t}}{\beta }%
\left( \mu _{t}^{\prime }+\mu _{t}\frac{h^{\prime }(t)}{h(t)}\right) .
\]%
Since the $\mu _{t}^{\prime }$ terms cancel on each side, this last holds iff%
\[
\psi _{t}\kappa =1-\frac{\psi _{t}}{\beta }\mu _{t}\frac{h^{\prime }(t)}{h(t)%
},\text{ so }-1+\psi _{t}\kappa =-\frac{\psi _{t}}{\beta }\mu _{t}\frac{%
h^{\prime }(t)}{h(t)}.
\]%
Computing the left-hand side, using the definition of $\psi _{t},$ gives%
\[
-1+\psi _{t}\kappa =-1+\frac{\kappa }{\kappa -\mu _{t}\lambda h(t)/\beta }=%
\frac{\mu _{t}\lambda h(t)}{\kappa -\mu _{t}\lambda h(t)/\beta }.
\]%
So, again using $\psi _{t},$
\[
\psi _{t}\mu _{t}\lambda h(t)=-\psi _{t}\mu _{t}\frac{h^{\prime }(t)}{h(t)},
\]%
implying, as $\psi _{t}>0$, either $\mu _{t}=0$ or $\lambda
h(t)^{2}=-h^{\prime }(t).$ The latter gives%
\[
\frac{dh}{h^{2}}=-\lambda dt:\text{ }h^{-1}=\lambda t+\text{const.}
\]%
But on $I$ this contradicts%
\[
h(t):=[2\Phi (\hat{\sigma}/2)-1],\text{ where }\hat{\sigma}=\sigma (1-t).
\]%
So consistency requires that $\mu _{t}=0$ on $I.$ But in this case the
co-state equation,%
\[
\pi _{t}[\mu _{t}\mathbb{\lambda }h(t)-(\beta -\alpha )]=\mu _{t}^{\prime
}-(\beta -\alpha ),
\]%
implies that $\pi _{t}=1$ on the interval $I$.\hfill $\square $

\section{Appendix B -- technicalities}

\subsection{Explicit valuation from the sparing policy}

Recall that%
\[
h(t):=[2\Phi (\hat{\sigma}/2)-1],\text{ where }\hat{\sigma}=\sigma (1-t),
\]%
so%
\[
\int_{0}^{s}\lambda \lbrack 2\Phi (\hat{\sigma}/2)-1]\text{ }\mathrm{d}%
t=2\lambda \int_{0}^{s}\Phi (\hat{\sigma}/2)\text{ }\mathrm{d}t-\lambda s.
\]

Put $a=\sigma /2$ and consider%
\begin{eqnarray*}
\int_{0}^{s}\Phi (a(1-t))\text{ }\mathrm{d}t &=&s\Phi
(a(1-s))+\int_{0}^{s}at\varphi (a(1-t))\text{ }\mathrm{d}t\text{ put }%
u=a(1-t)=a-at \\
&=&s\Phi (a(1-s))+\frac{1}{a}\int_{a}^{a(1-s)}(u-a)\varphi (u)\text{ }%
\mathrm{d}u\text{ as }\mathrm{d}u=-a\text{ }\mathrm{d}t \\
&=&s\Phi (a(1-s))+\frac{1}{a}\int_{a}^{a(1-s)}u\varphi (u)\text{ }\mathrm{d}%
u-[\Phi (a(1-s))-\Phi (a)] \\
&=&\Phi (a)+(s-1)\Phi (a(1-s))+\frac{1}{a}\int_{a}^{a(1-s)}u\varphi (u)\text{
}\mathrm{d}u\text{ (see below)} \\
&=&\Phi (a)+(s-1)\Phi (a(1-s))+\frac{1}{a\sqrt{2\pi }}%
[-e^{-u^{2}/2}]_{a}^{a(1-s)} \\
&=&\Phi (a)+(s-1)\Phi (a(1-s))+\frac{1}{a\sqrt{2\pi }}%
[e^{-a^{2}/2}-e^{-a^{2}(1-s)^{2}/2}].
\end{eqnarray*}%
Put%
\[
g(s):=\int_{0}^{s}h(t)\text{ }\mathrm{d}t=\int_{0}^{s}[2\Phi (\hat{\sigma}%
/2)-1]\text{ }\mathrm{d}t
\]%
\[
=2\left[ \Phi (a)+(s-1)\Phi (a(1-s))+\frac{1}{a\sqrt{2\pi }}%
[e^{-a^{2}/2}-e^{-a^{2}(1-s)^{2}/2}]\right] -s.
\]%
Thus $g(0)=0$ and $g(1)=2\left[ \Phi (a)+\frac{1}{a\sqrt{2\pi }}%
[e^{-a^{2}/2}-1]\right] -1\simeq \frac{a}{\sqrt{2\pi }}.$

For small $a,$%
\begin{eqnarray*}
g(t) &\simeq &2\left[ \frac{1}{2}+(t-1)\frac{1}{2}+\frac{1}{a\sqrt{2\pi }}%
[e^{-a^{2}/2}-e^{-a^{2}(1-s)^{2}/2}]\right] -t \\
&\simeq &t+\frac{a}{\sqrt{2\pi }}.
\end{eqnarray*}%
So
\begin{eqnarray*}
\gamma _{t}^{1} &=&\exp (-\lambda g(t)) \\
&=&\exp \left( -2\lambda \left[ \Phi (a)+(t-1)\Phi (a(1-t))+\frac{1}{a\sqrt{%
2\pi }}[e^{-a^{2}/2}-e^{-a^{2}(1-t)^{2}/2}]\right] -\lambda t\right) .
\end{eqnarray*}

\subsection{Cascading algorithm}

Here we characterize optimal switching points as solutions of equations that
depend only on earlier switching points. To do this we must first describe
the very simple form of the trajectory $\gamma _{t}$ when its behaviour is
controlled by an alternating policy, switching consecutively between $\pi =0$
and $\pi =1.$ This expands the Remark in Section 3.3 on the behaviour of $%
\gamma _{t}.$ We use the notation $\gamma _{i}=\gamma _{\theta _{i}}$ for
the value of $\gamma _{t}$ at the $i$-th switching point $\theta _{i};$
notice that if%
\[
\gamma _{i}=\gamma _{i+1},
\]%
as when $\pi =0$ on $[\theta _{i},\theta _{i+1}),$ then $\gamma _{i+2}$
acquires two exponential factors:%
\[
\gamma _{i+2}=\gamma _{i+1}e^{\lambda g(\theta _{i+1})}e^{-\lambda g(\theta
_{i+2})}<\gamma _{i+1}.
\]

\bigskip

\noindent \textbf{Lemma 4 (Behaviour of }$\gamma _{t}$\textbf{). }\textit{%
Set }$\theta _{0}=0$\textit{\ and }$\theta _{n+1}=1,$\textit{\ and assume
consecutive switches occur at time locations }$\theta _{1}<\theta
_{2}<...<\theta _{n}<1.$

\textit{\noindent }(i)\textit{\ If }$\pi =0$\textit{\ near }$t=0,$\textit{\
so that }$1=\gamma _{0}=\gamma _{\theta _{1}},$\textit{\ then for even }$%
i\leq n$\textit{\ }%
\[
\gamma _{\theta _{i}}=\gamma _{\theta _{i+1}}=e^{\lambda g(\theta
_{1})}e^{-\lambda g(\theta _{2})}e^{\lambda g(\theta _{3})}e^{-\lambda
g(\theta _{4})}...e^{\lambda g(\theta _{i-1})}e^{-\lambda g(\theta _{i})},
\]%
\textit{e.g.}%
\[
\gamma _{\theta _{2}}=e^{\lambda g(\theta _{1})}e^{-\lambda g(\theta _{2})}%
\text{ as }\gamma _{t}=e^{-\lambda g(t)}e^{\lambda g(\theta _{1})}\text{ for
}t\in \lbrack \theta _{1},\theta _{2}].
\]

\textit{\noindent }(ii)\textit{\ If }$\pi =1$\textit{\ near }$t=0,$\textit{\
so that }$\gamma _{\theta _{2}}=\gamma _{\theta _{1}}=e^{-\lambda g(\theta
_{1})},$\textit{\ then for odd }$i\leq n$\textit{\ }%
\[
\gamma _{\theta _{i+1}}=\gamma _{\theta _{i}}=e^{-\lambda g(\theta
_{1})}e^{\lambda g(\theta _{2})}e^{-\lambda g(\theta _{3})}e^{\lambda
g(\theta _{4})}...e^{\lambda g(\theta _{i-1})}e^{-\lambda g(\theta _{i})},
\]%
\textit{e.g.}%
\[
\gamma _{\theta _{3}}=e^{-\lambda g(\theta _{1})}e^{\lambda g(\theta
_{2})}e^{-\lambda g(\theta _{3})}\text{ as }\gamma _{t}=\gamma _{\theta
_{2}}e^{\lambda g(\theta _{2})}e^{-\lambda g(t)}\text{ for }t\in \lbrack
\theta _{2},\theta _{3}].
\]

\noindent \textbf{Proof. }Suppose consecutive switches occur at some time
locations $\theta _{1}^{\prime }<\theta _{2}^{\prime }<\theta _{3}^{\prime }$
with $\pi =1$ on $(\theta _{1}^{\prime },\theta _{2}^{\prime })$ and $\pi =0$
on $(\theta _{2}^{\prime },\theta _{3}^{\prime }).$ Then $\gamma
_{t}=e^{\lambda g(\theta _{1}^{\prime })}e^{-\lambda g(t)}\gamma _{\theta
_{1}^{\prime }}$ for $t\in \lbrack \theta _{1}^{\prime },\theta _{2}^{\prime
}],$ so that
\[
\gamma _{\theta _{3}^{\prime }}=\gamma _{\theta _{2}^{\prime }}=\gamma
_{\theta _{1}^{\prime }}e^{\lambda g(\theta _{1}^{\prime })}e^{-\lambda
g(\theta _{2}^{\prime })}.
\]%
If, however, $\pi =0$ on $(\theta _{1}^{\prime },\theta _{2}^{\prime })$ and
$\pi =1$ on $(\theta _{2}^{\prime },\theta _{3}^{\prime }),$ then $\gamma
_{\theta _{1}^{\prime }}=\gamma _{\theta _{2}^{\prime }}$ and $\gamma
_{t}=\gamma _{\theta _{2}^{\prime }}e^{-\lambda g(t)}e^{\lambda g(\theta
_{2}^{\prime })}$ for $t\in $ $[\theta _{2}^{\prime },\theta _{3}^{\prime
}], $ so%
\[
\gamma _{\theta _{3}^{\prime }}=\gamma _{\theta _{2}^{\prime }}e^{\lambda
g(\theta _{2}^{\prime })}e^{-\lambda g(\theta _{3}^{\prime })}=\gamma
_{\theta _{1}^{\prime }}e^{\lambda g(\theta _{2}^{\prime })}e^{-\lambda
g(\theta _{3}^{\prime })}.
\]%
Both formulas have alternating signs with the negatively signed exponential
corresponding to the endpoint of $\pi =1$ constancy; both relate $\gamma
_{\theta _{3}^{\prime }}$ to $\gamma _{\theta _{1}^{\prime }},$ i.e. to two
switches earlier.

To generalizing to (exactly) $n$ consecutive switches occuring at times $%
\theta _{1}<\theta _{2}<...<\theta _{n}<1,$ apply the preceding formula
inductively backwards.

If $\pi =0$ near $t=0,$ so that $1=\gamma _{0}=\gamma _{\theta _{1}},$ then
with $i\leq n$ \textit{even, }$\gamma _{t}$ is constant on $[\theta
_{i},\theta _{i+1}].$ So $\pi =1$ on $(\theta _{1}^{\prime },\theta
_{2}^{\prime })=(\theta _{i-1},\theta _{i})$ and $\pi =0$ on $(\theta
_{2}^{\prime },\theta _{3}^{\prime })=(\theta _{i},\theta _{i+1}).$ With due
regard to parity, we obtain, since $\theta _{i}$ is an endpoint for $\pi =1$
constancy,%
\[
\gamma _{\theta _{i+1}}=\gamma _{\theta _{i}}=\gamma _{\theta
_{i-1}}e^{\lambda g(\theta _{i-1})}e^{-\lambda g(\theta _{i})}.
\]%
Furthermore, $\pi =0$ on $(\theta _{i-2},\theta _{i-1})$ so%
\[
\gamma _{\theta _{i}}=\gamma _{\theta _{i-1}}e^{\lambda g(\theta
_{i-1})}e^{-\lambda g(\theta _{i})}=\gamma _{\theta _{i-2}}e^{\lambda
g(\theta _{i-1})}e^{-\lambda g(\theta _{i})}
\]
and $i-2$ is even. Hence%
\[
\gamma _{\theta _{i+1}}=\gamma _{\theta _{i}}=e^{\lambda g(\theta
_{1})}e^{-\lambda g(\theta _{2})}e^{\lambda g(\theta _{3})}e^{-\lambda
g(\theta _{4})}...e^{\lambda g(\theta _{i-1})}e^{-\lambda g(\theta _{i})},
\]%
valid because for $i=0$ this yields $\gamma _{\theta _{0}}=1,$ since $%
g(\theta _{i})=g(0)=0.$

Likewise, if $\pi =1$ near $t=0,$ so that $\gamma _{\theta _{2}}=\gamma
_{\theta _{1}}=e^{-\lambda g(\theta _{1})},$ then with $i\leq n$ \textit{odd
}$\gamma _{t}$ is constant on $[\theta _{i},\theta _{i+1}]$ so that $\pi =0$
on $(\theta _{i},\theta _{i+1})$ and so $\pi =1$ on $(\theta _{i-1},\theta
_{i}).$ Thus
\[
\gamma _{\theta _{i}}=\gamma _{\theta _{i-1}}e^{\lambda g(\theta
_{i-1})}e^{-\lambda g(\theta _{i})}=\gamma _{\theta _{i-2}}e^{\lambda
g(\theta _{i-1})}e^{-\lambda g(\theta _{i})}
\]
as $\pi =0$ on $(\theta _{i-2},\theta _{i-1})$ and here $i-2$ is odd. Hence
\[
\gamma _{\theta _{i+1}}=\gamma _{\theta _{i}}=e^{-\lambda g(\theta
_{1})}e^{\lambda g(\theta _{2})}e^{-\lambda g(\theta _{3})}e^{\lambda
g(\theta _{4})}...e^{\lambda g(\theta _{i-1})}e^{-\lambda g(\theta _{i})}.
\]%
This is valid since for $i=1$ the formula reduces to $\gamma _{\theta
_{i}}=e^{-\lambda g(\theta _{1})}$.\hfill $\square $

\bigskip

\noindent \textbf{Remark.} If $\pi =0$ near $t=0$ we thus have%
\begin{eqnarray*}
\gamma _{\theta _{1}} &=&\gamma _{0}=1,\gamma _{\theta _{3}}=\gamma _{\theta
_{2}}=e^{\lambda \lbrack g(\theta _{1})-g(\theta _{2})]}, \\
\gamma _{\theta _{5}} &=&\gamma _{\theta _{4}}=e^{\lambda \lbrack g(\theta
_{1})-g(\theta _{2})+g(\theta _{3})-g(\theta _{4})]},...,
\end{eqnarray*}%
so that%
\[
\log \gamma _{2n}=\lambda \tsum\nolimits_{i=0}^{2n}(-1)^{i+1}g(\theta _{i})%
\text{ with }\theta _{0}=0.
\]%
Similarly, if $\pi =1$ near $t=0,$ we have
\[
\gamma _{0}=1,\text{ }\gamma _{\theta _{2}}=\gamma _{\theta
_{1}}=e^{-\lambda g(\theta _{1})},\text{ }\gamma _{\theta _{4}}=\gamma
_{\theta _{3}}=e^{-\lambda g(\theta _{1})}e^{\lambda g(\theta
_{2})}e^{-\lambda g(\theta _{3})},...,
\]%
so that%
\[
\log \gamma _{2n+1}=\lambda \tsum\nolimits_{i=0}^{2n+1}(-1)^{i}g(\theta _{i})%
\text{ with }\theta _{0}=0.
\]%
In either case $\gamma _{i}$ depends either on $\{\theta _{j}:j\leq i\}$ or
on $\{\theta _{j}:j<i\}.$ More specifically, if $\pi =0$ near $t=0$ then $%
\gamma _{i}$ for even $i$ depends on $\{\theta _{j}:j\leq i\}.$ Similarly,
if $\pi =1$ near $t=0,$ then then $\gamma _{i}$ for odd $i$ depends on $%
\{\theta _{j}:j\leq i\}.$

In both cases these $\gamma $ values are monotonically decreasing since%
\[
g(\theta _{i})<g(\theta _{i+1})\text{ so that }e^{\lambda \lbrack g(\theta
_{i})-g(\theta _{i+1})]}<1
\]
with $\gamma _{2j}$ following $\pi =0$ (respectively $\gamma _{2j+1}$
following $\pi =1)$ being strictly monotone.

\bigskip

\noindent \textbf{Definitions.} Put%
\[
\gamma _{i}=\gamma _{\theta _{i}}\text{ with }\gamma _{0}=1.
\]%
Furthermore, put $\kappa _{0}=\kappa $ and for $i\geq 1$%
\[
\kappa _{i}=\kappa \gamma _{i}.
\]%
Thus if $\pi =0$ on $[0,\theta _{1}],$ then $\gamma _{0}=\gamma
_{1}=1,\gamma _{2}=e^{\lambda \lbrack g(\theta _{1})-g(\theta _{2})]}$ and
so the $\kappa $-sequence is%
\[
\kappa _{0}=\kappa ,\text{ }\kappa _{1}=\kappa ,\text{ }\kappa
_{2}=e^{-\lambda \lbrack g(\theta _{2})-g(\theta _{1})]}\kappa ,...
\]%
and if $\pi =1$ on $[0,\theta _{1}],$ then $\gamma _{0}=1$ and $\gamma
_{2}=_{0}\gamma _{1}=e^{-\lambda g(\theta _{1})}$ and so the $\kappa $%
-sequence is%
\begin{eqnarray*}
\kappa _{0} &=&\kappa ,\text{ }\kappa _{1}=\kappa e^{-\lambda g(\theta
_{1})};\text{ }\kappa _{2}=_{0}\kappa _{1}=e^{-\lambda g(\theta _{1})}; \\
\kappa _{3} &=&e^{-\lambda g(\theta _{1})}e^{-\lambda \lbrack g(\theta
_{3})-g(\theta _{2})]};\kappa _{4}=_{0}\kappa _{3};...
\end{eqnarray*}

Note. Either $\kappa _{i}=_{\pi =0}\kappa _{i-1}$ or $\kappa _{i}=_{\pi
=1}\kappa _{i-1}e^{-\lambda \lbrack g(\theta _{i})-g(\theta _{i-1})]}$ using
the factor $e^{-\lambda \lbrack g(t)-g(\theta _{i-1})]}$

\bigskip

\bigskip

\noindent \textbf{Proof of Theorem 3} We consider the
impact of $\theta =\theta _{i}$ on the objective function, as measured by
\[
I_{i}(\pi ,\theta )/\beta =\int_{\theta _{i-1}}^{\theta }(1-\pi _{t})[\gamma
_{t}^{1}-\kappa \gamma _{t}]\text{ }\mathrm{d}t+\int_{\theta }^{\theta
_{i+1}}(1-\pi _{t})[\gamma _{t}^{1}-\kappa \gamma _{t}]\text{ }\mathrm{d}t.
\]%
Since $\theta _{0}=0$ and $\theta _{n+1}=1,$%
\[
I(\pi )/\beta =\tsum\nolimits_{i=1}^{n}I_{i}(\pi ,\theta _{i})/\beta .
\]

If $\pi =0$ on $[\theta _{i-1},\theta _{i}],$then $\gamma _{t}\equiv \gamma
_{\theta _{i-1}}=\gamma _{i-1}$ and $\gamma _{t}^{1}=e^{-\lambda g(t)}$ on
this interval and so the impact on the objective function reduces to%
\[
I_{i}(\pi ,\theta )/\beta =\int_{\theta _{i-1}}^{\theta }e^{-\lambda g(t)}%
\text{ }\mathrm{d}t-\gamma _{i-1}\kappa (\theta -\theta _{i-1}),
\]%
since $\pi =1$ on $[\theta _{i},\theta _{i+1}]$ yields a zero contribution.
Differentiation w.r.t. $\theta $ yields the first-order condition%
\[
e^{-\lambda g(\theta )}-\gamma _{i-1}\kappa =0,\text{ i.e. }e^{-\lambda
g(\theta _{i})}=\gamma _{i}\kappa ,
\]%
since $\gamma _{i}=\gamma _{i-1},$ as $\gamma $ is constant on $[\theta
_{i-1},\theta _{i}]$.

(ii) If $\pi =1$ on $[\theta _{i-1},\theta _{i}]$, then $\gamma
_{t}^{1}\equiv e^{-\lambda g(t)}$ as always and $\gamma _{t}=\gamma _{\theta
_{i-1}}e^{\lambda \lbrack g(\theta _{i-1})-g(t)]}$ on this interval. So
taking $\theta _{i}=\theta $ with $\theta _{i-1}<\theta <\theta _{i+1}$ we
have%
\[
\gamma _{s}=\gamma _{\theta _{i-1}}e^{\lambda g(\theta _{i-1})}e^{-\lambda
g(s)}\text{ for }\theta _{i-1}<s<\theta _{i},\text{ as }\pi =1\text{ thereon.%
}
\]%
So for $t\in \lbrack \theta ,\theta _{i+1})$
\begin{eqnarray*}
\gamma _{t} &\equiv &\gamma _{\theta _{i-1}}e^{\lambda g(\theta
_{i-1})}e^{-\lambda g(\theta )}\text{ as }\pi =0\text{ on }[\theta ,\theta
_{i+1}),\text{ since } \\
\gamma _{s} &=&\gamma _{\theta _{i-1}}e^{\lambda g(\theta
_{i-1})}e^{-\lambda g(s)}\text{ for }s\in \lbrack \theta _{i-1},\theta )%
\text{ (where }\pi =1\text{).}
\end{eqnarray*}%
To abbreviate the (constant) expressions independent of $\theta $, we will
write%
\[
\gamma _{t}\equiv \tilde{\gamma}e^{-\lambda g(\theta )}\text{ with }\pi (t)=0%
\text{ for }t\in \lbrack \theta ,\theta _{i+1}).
\]%
Here the impact on the objective function reduces to%
\begin{eqnarray*}
I_{i}(\pi ,\theta )/\beta &=&\int_{\theta }^{\theta _{i+1}}[\gamma
_{t}^{1}-\kappa \tilde{\gamma}e^{-\lambda g(\theta )}]\text{ }\mathrm{d}t \\
&=&\int_{\theta }^{\theta _{i+1}}e^{-\lambda g(t)}\text{ }\mathrm{d}t-\kappa
\tilde{\gamma}e^{-\lambda g(\theta )}(\theta _{i+1}-\theta ),
\end{eqnarray*}%
since there is a zero contribution to the objective function on $[\theta
_{i-1},\theta ],$ where $\pi =1$. Differentiation w.r.t. $\theta $ yields
the first-order condition at $\theta =\theta _{i}$%
\[
-e^{-\lambda g(\theta )}+\kappa \tilde{\gamma}e^{-\lambda g(\theta )}+\kappa
\tilde{\gamma}e^{-\lambda g(\theta )}\lambda h(\theta )(\theta _{i+1}-\theta
)=0,
\]%
or%
\begin{eqnarray*}
\kappa \tilde{\gamma}e^{-\lambda g(\theta )}\lambda h(\theta )(\theta
_{i+1}-\theta ) &=&e^{-\lambda g(\theta )}-\kappa \tilde{\gamma}e^{-\lambda
g(\theta )}, \\
\kappa \tilde{\gamma}\lambda h(\theta )(\theta _{i+1}-\theta ) &=&1-\kappa
\tilde{\gamma},
\end{eqnarray*}%
so%
\[
\theta _{i+1}=\theta _{i}+\frac{1-\kappa \tilde{\gamma}}{\kappa \tilde{\gamma%
}\lambda h(\theta _{i})}.
\]%
Now $\kappa _{i}=\kappa \gamma _{i}=\kappa \gamma _{\theta _{i}}=\kappa
\tilde{\gamma}e^{-\lambda g(\theta _{i})}$ (as above), so $\kappa \tilde{%
\gamma}=\kappa _{i}e^{\lambda g(\theta _{i})},$ so a re-arrangement yields
the claim.\hfill $\square $

\bigskip

\noindent \textbf{Corollary 5. }\textit{For each switching point }$\theta
_{i}$\textit{\ that is a right endpoint of an interval with }$\pi =1$\textit{%
\ constancy, i.e. with }$\pi =1$ \textit{on }$[\theta _{i-1},\theta _{i})$%
\textit{,}%
\[
g\left( \theta _{i}+\frac{1-\kappa _{i-1}e^{\lambda g(\theta _{i-1})}}{%
\kappa _{i-1}e^{\lambda g(\theta _{i-1})}\lambda h(\theta _{i})}\right)
=g(\theta _{i})-(\log \kappa _{i-1})/\lambda -g(\theta _{i-1}).
\]%
\textit{where }$\kappa _{i-1}$\textit{\ depends on }$\{\theta _{j}:j\leq
i-1\}.$\textit{\ Its solution defines}%
\[
\theta _{i+1}=\theta _{i}+\frac{1-\kappa _{i-1}e^{\lambda g(\theta _{i-1})}}{%
\kappa _{i-1}e^{\lambda g(\theta _{i-1})}\lambda h(\theta _{i})},
\]%
\textit{If }$\theta _{1}$\textit{\ is a right endpoint of an interval with }$%
\pi =0$\textit{\ constancy, then}%
\[
g(\theta _{1})=-\log \kappa /\lambda .
\]%
\textit{In particular, if }$\pi =0$\textit{\ on }$[0,\theta _{1}),$\textit{\
so that }$\gamma _{1}=\gamma _{0}=1,$\textit{\ then }$\pi =1$\textit{\ on }$%
[\theta _{1},\theta _{2})$\textit{\ and so for }$i=2,$\textit{\ as }$\kappa
_{1}=\kappa =e^{-\lambda g(\theta _{1})},$\textit{\ }%
\[
\theta _{3}=\theta _{2}+\frac{1-\kappa _{1}e^{\lambda g(\theta _{1})}}{%
\kappa _{1}e^{\lambda g(\theta _{1})}\lambda h(\theta _{2})}=\theta _{2},
\]%
\textit{a contradiction. Consequently, there cannot be a further switching
point.}

\textit{Furthermore, the sequence }$\kappa _{i}$\textit{\ is (weakly)
decreasing with alternate members strictly decreasing.}

\bigskip

\noindent \textbf{Proof of Corollary 5. }Consider a pair of intervals: $%
[\theta _{i-1},\theta _{i})$ on which $\pi =1,$ contiguous with $[\theta
_{i},\theta _{i+1})$ on which $\pi =0.$ By Lemma 4, $\gamma _{i}=\gamma
_{i-1}e^{\lambda g(\theta _{i-1})-\lambda g(\theta _{i})}$. Since $\pi =1$
on $(\theta _{i-1},\theta _{i})$ and%
\begin{eqnarray*}
\kappa _{i}e^{\lambda g(\theta _{i})} &=&\kappa \gamma _{i}e^{\lambda
g(\theta _{i})}=\kappa \gamma _{i-1}e^{\lambda g(\theta _{i-1})-\lambda
g(\theta _{i})}e^{\lambda g(\theta _{i})}=\kappa \gamma _{i-1}e^{\lambda
g(\theta _{i-1})} \\
&=&\kappa _{i-1}e^{\lambda g(\theta _{i-1})}.
\end{eqnarray*}%
By Theorem 2$^{\prime }$%
\[
\theta _{i+1}=\theta _{i}+\frac{1-\kappa _{i}e^{\lambda g(\theta _{i})}}{%
\kappa _{i}\lambda h(\theta _{i})e^{\lambda g(\theta _{i})}}=\theta _{i}+%
\frac{1-\kappa _{i-1}e^{\lambda g(\theta _{i-1})}}{\kappa _{i-1}e^{\lambda
g(\theta _{i-1})}\lambda h(\theta _{i})},
\]%
or equivalently%
\[
\kappa _{i-1}=\frac{e^{-\lambda g(\theta _{i-1})}}{1+\lambda (\theta
_{i+1}-\theta _{i})h(\theta _{i})}.
\]%
Also $\kappa _{i+1}=\kappa _{i}=\kappa \gamma _{i}$ as $\pi =0$ on $[\theta
_{i},\theta _{i+1}),$ so by Theorem 2$^{\prime }$ (and since $\gamma
_{i}=\gamma _{i-1}e^{\lambda g(\theta _{i-1})-\lambda g(\theta _{i})},$ as
noted earlier)
\begin{eqnarray*}
\text{ }e^{-\lambda g(\theta _{i+1})} &=&\kappa _{i+1}=\kappa _{i}=\kappa
\gamma _{i}=\kappa _{i-1}e^{-\lambda \lbrack g(\theta _{i})-g(\theta
_{i-1})]}:\text{ } \\
-\lambda g(\theta _{i+1}) &=&\log \kappa _{i-1}-\lambda g(\theta
_{i})+\lambda g(\theta _{i-1}).
\end{eqnarray*}%
From here%
\[
g(\theta _{i+1})=-\frac{\log \kappa _{i-1}}{\lambda }+g(\theta
_{i})-g(\theta _{i-1})\text{ with }\pi =0\text{ on }[\theta _{i},\theta
_{i+1}).
\]%
But at the same time%
\[
g(\theta _{i-1})=-\log (\kappa _{i-1})/\lambda ,\text{ as also }\pi =0\text{
on }[\theta _{i-2},\theta _{i-1}),
\]%
so%
\[
g(\theta _{i+1})=-\frac{\log \kappa _{i-1}}{\lambda }+g(\theta
_{i})-g(\theta _{i-1})\text{ with }\pi =0\text{ on }[\theta _{i},\theta
_{i+1}).
\]%
\[
g(\theta _{i+1})=g(\theta _{i})
\]

This equation coupled with the earlier one yields the claim.\hfill $\square $

\bigskip

\noindent \textbf{Remark.} Other than the single equation in one variable
when the equilibrium demands that $\pi =0$ is applied on $[0,\theta _{1}),$
the optimization conditions yield linked pairs of simultaneous equations
reduceable to a single equation involving only earlier switching points.

\bigskip

\textbf{Remark on piecewise constant }$\kappa _{t}$\textbf{.} The argument
proving Theorem 2$^{\prime }$ about optimal switching points, allows $\kappa
$ to be piecewise constant with constancy between switching points. An
example is illustrated in Figure 3 in Section 4.2. For suppose we have $%
\kappa =\kappa _{-}=\kappa _{\theta _{i-1}}$ to the left of $\theta _{i}$
and $\kappa =\kappa _{+}=\kappa _{\theta _{i}}$ to the right (i.e. c\`{a}dl%
\`{a}g as before) that $\kappa =\kappa _{\theta _{i}}$ on $[\theta
_{i},\theta _{i+1})$), then the corresponding formulas remain the same with
only a change to the definition of $\kappa _{i},$ thus%
\[
e^{-\lambda g(\theta _{i})}=\bar{\kappa}_{i}:=\gamma _{i}\kappa _{-}=\gamma
_{i}\kappa _{\theta _{i-1}}\text{ assuming }\pi =0\text{ on }[\theta
_{i-1},\theta _{i})
\]%
and%
\begin{eqnarray*}
\theta _{i+1} &=&\theta _{i}+\frac{1-\kappa _{+}\gamma _{i}e^{\lambda
g(\theta _{i})}}{\kappa _{+}\gamma _{i}e^{\lambda g(\theta _{i})}\lambda
h(\theta _{i})} \\
&=&\theta _{i}+\frac{\bar{\kappa}_{i}^{-1}e^{-\lambda g(\theta _{i})}-1}{%
\lambda h(\theta _{i})}\text{ with }\bar{\kappa}_{i}:=\gamma _{i}\kappa
_{+}=\gamma _{i}\kappa _{\theta _{i}},
\end{eqnarray*}%
assuming $\pi =1$ on $[\theta _{i-1},\theta _{i})$, so that $\gamma
_{i}=\gamma _{\theta _{i-1}}e^{\lambda \lbrack g(\theta _{i-1})-g(\theta
_{i})]}.$

However, monotonicity claims concerning $\bar{\kappa}_{i}$ would now rely on
assumptions about the sequence $\kappa _{\theta _{i}}.$

\bigskip

Indeed, suppose $\pi =1$ is intended for $t>\theta _{i},$ then we can
interpret $\kappa $ as being $\kappa _{-}$ to the right of $\theta _{i}$
without altering the payoff (because of the factor $1-\pi _{t}).$ We then get%
\[
e^{-\lambda g(\theta _{i})}=\gamma _{i}\kappa _{-}.
\]%
Similarly, if $\pi =0$ is intended for $t>\theta _{i},$ i.e. $\pi =1$ is
intended for $t<\theta _{i},$ then we can interpret $\kappa $ as being $%
\kappa _{+}$ to the left of $\theta _{i}$ again without altering the payoff.
We then get%
\[
\theta _{i+1}=\theta _{i}+\frac{1-\kappa _{+}\tilde{\gamma}}{\kappa _{+}%
\tilde{\gamma}\lambda h(\theta _{i})}.
\]%
The behaviour of $\gamma _{t}$ near $\theta _{i}$ is determined by the ODE
and by the value of $\theta _{i}$. So since $\tilde{\gamma}=\gamma _{\theta
_{i-1}}e^{\lambda g(\theta _{i-1})}=\gamma _{i}e^{\lambda g(\theta _{i})}$
as $\gamma _{\theta _{i}}=\gamma _{\theta _{i-1}}e^{\lambda \lbrack g(\theta
_{i-1})-g(\theta _{i})]}$ one arrives at%
\[
\theta _{i+1}=\theta _{i}+\frac{1-\kappa _{+}\gamma _{i}e^{\lambda g(\theta
_{i})}}{\kappa _{+}\gamma _{i}e^{\lambda g(\theta _{i})}\lambda h(\theta
_{i})}.
\]%
This agrees for $\kappa $ constant with the previous context.

\bigskip

Finally, we show that there are critical pairs $\kappa ,\lambda $ with
insoluble paired linking equations. Notice that for $0<\eta <h(0)$ the map%
\[
\lambda \rightarrow \kappa (\lambda )=(1+\eta \lambda )^{-1}
\]%
is decreasing in $\lambda $ and also in $\eta $. For large enough $\lambda ,$
this map will assigns to $\lambda $ a small enough kappa for which the
linking equation must fail. We may thus term%
\[
\lambda \rightarrow \kappa (\lambda )=(1+\lambda h(0))^{-1}
\]%
the \textit{killing map} for $\lambda .$ At that level of $\kappa $ there
cannot be any further switching. The proof uses compactness (equivalently
uniform continuity).

\bigskip

\noindent \textbf{Lemma 5 (Lemma on killing).} \textit{Take }$0<\eta <h(0)$
\textit{and}%
\[
T(\theta )=\theta +\eta /h(\theta ),
\]%
\textit{which is strictly increasing and continuous, and put }$\bar{\theta}%
_{1}=T^{-1}(1).$\textit{\ }

\textit{Then there exists }$\kappa _{c}$\textit{\ such that for all }$\kappa
\leq \kappa _{c}$\textit{\ and for }$\lambda _{\kappa }:=(\kappa
^{-1}-1)/\eta \geq \lambda _{c}:=(\kappa _{c}^{-1}-1)/\eta $%
\[
g(\theta )+\frac{\log \kappa ^{-1}}{\lambda _{\kappa }}<g\left( \theta +%
\frac{(\kappa ^{-1}-1)}{\lambda _{\kappa }h(\theta )}\right) \text{ for all }%
\theta \in \lbrack 0,\bar{\theta}_{1}].
\]%
\textit{Equivalently there exists }$\lambda _{c}$\textit{\ such that if }%
\[
\kappa (\lambda ):=(1+\eta \lambda )^{-1}\leq \kappa (\lambda _{c}),
\]%
\textit{then with }$\kappa =\kappa (\lambda )$%
\[
g(\theta )+\frac{\log \kappa ^{-1}}{\lambda }<g\left( \theta +\frac{(\kappa
^{-1}-1)}{\lambda h(\theta )}\right) \text{ for all }\theta \in \lbrack 0,%
\bar{\theta}_{1}].
\]

\bigskip

\noindent \textbf{Proof.} $T$ is continuous and strictly increasing, hence
so is $T^{-1}.$ Since%
\[
\theta +\frac{\eta }{h(\theta )}=1\text{ iff }\eta =(1-\theta )h(\theta ),
\]%
so $T(\theta )=1$ is soluble with unique solution $\theta =\bar{\theta}_{1}$
as $(1-\theta )h(\theta )$ decreases from $h(0)$ to $0$ at $\theta =1.$Thus $%
T$ takes $0$ to $\eta /h(0)$ and $\bar{\theta}_{1}$ to $1.$ As $T$ is
strictly increasing, $T^{-1}$ is well defined and takes $[\eta /h(0),1]$ to $%
[0,\bar{\theta}_{1}].$

For $\eta /h(0)\leq k\leq 1$ let $\theta =\bar{\theta}_{k}$ solve $T(\theta
)=k,$ i.e. $\bar{\theta}_{k}=T^{-1}(k):$%
\[
T(\theta )=\theta +\frac{\eta }{h(\theta )}=k.
\]

Clearly $\bar{\theta}_{k}<k.$ So, equivalently, $\theta =\bar{\theta}_{k}<k$
and solves in $[0,\bar{\theta}_{1}]$
\[
\eta =(k-\theta )h(\theta ),
\]%
where the function on the right is decreasing for $\theta <k$, ranging from $%
kh(0)\geq \eta $ down to $0,$ so $\bar{\theta}_{k}=\bar{\theta}_{k}(\eta )$
is well-defined and the range of $k\rightarrow $ $\bar{\theta}_{k}$ on $%
[\eta /h(0),1]$ is $[0,\bar{\theta}_{1}]$

Note that $\bar{\theta}_{\eta /h(0)}=0$. Now, for $1\geq k\geq \eta /h(0),$%
\[
g(\bar{\theta}_{k})<g(T(\bar{\theta}_{k}))=g(k)\leq g(1)
\]%
(as $g$ is strictly increasing), equivalently%
\[
g(T^{-1}(k))<g(k).
\]%
So we consider the stronger condition%
\[
g(\bar{\theta}_{k})+\eta \frac{\log \kappa ^{-1}}{(\kappa ^{-1}-1)}<g(k)=g(T(%
\bar{\theta}_{k}(\eta ))),
\]%
which holds provided we choose $\kappa =\kappa _{k}\in (0,1)$ with%
\[
\frac{\log \kappa _{k}^{-1}}{(\kappa _{k}^{-1}-1)}<[g(k)-g(\bar{\theta}%
_{k})]/\eta
\]%
and this holds for all $\kappa <\kappa _{k}$, since
\[
\lim_{\kappa \rightarrow 0}\frac{-\log \kappa }{\kappa ^{-1}-1}=\lim_{\kappa
\rightarrow 0}\frac{-\kappa ^{-1}}{-\kappa ^{-2}}=0,\text{ and }\lim_{\kappa
\rightarrow \infty }\frac{-\log \kappa }{\kappa ^{-1}-1}=\infty ,
\]%
and $\kappa \rightarrow -\log \kappa /(\kappa ^{-1}-1)$ is a strictly
increasing function.

Thus there is $\kappa _{k}$ such that
\[
g(T^{-1}(k))+\eta \frac{\log \kappa _{k}^{-1}}{(\kappa _{k}^{-1}-1)}<g(k).
\]

By continuity at $k$ of $g(T^{-1}(.))$ and of $g()$, for some $\delta
_{k}>0, $ there is an interval $J_{k}:=(k-\delta _{k},k+\delta _{k})\cap
\lbrack \eta /h(0),1]$ such that%
\[
g(T^{-1}(s))+\eta \frac{\log \kappa _{k}^{-1}}{(\kappa _{k}^{-1}-1)}<g(s)%
\text{ for all }s\in J_{k}.
\]%
By compactness of $[\eta /h(0),1]$, there is a finite set $K\subset \lbrack
\eta /h(0),1]$ such that the finite number of intervals $\{J_{k}:k\in K\}$
covers $[\eta /h(0),1].$

Let $\kappa _{c}=\min \{\kappa _{k}:k\in K\}.$ Then
\[
g(T^{-1}(s))+\eta \frac{\log \kappa _{c}^{-1}}{(\kappa _{c}^{-1}-1)}<g(s)%
\text{ for }s\in \lbrack \eta /h(0),1].
\]%
or writing $s=T(\theta )$ for $\theta \in \lbrack 0,\bar{\theta}_{1}],$ we
have%
\[
g(\theta )+\eta \frac{\log \kappa _{c}^{-1}}{(\kappa _{c}^{-1}-1)}%
<g(T(\theta ))\text{ for }\theta \in \lbrack 0,\bar{\theta}_{1}].
\]%
So for all $\kappa \leq \kappa _{c}$%
\[
g(\theta )+\eta \frac{\log \kappa ^{-1}}{(\kappa ^{-1}-1)}<g\left( \theta +%
\frac{\eta }{h(\theta )}\right) \text{ for }\theta \in \lbrack 0,\bar{\theta}%
_{1}].
\]%
For any such $\kappa $ take%
\[
\lambda _{\kappa }:=(\kappa ^{-1}-1)/\eta \geq (\kappa _{c}^{-1}-1)/\eta .
\]%
then%
\[
g(\theta )+\frac{\log \kappa ^{-1}}{\lambda _{\kappa }}<g\left( \theta +%
\frac{(\kappa ^{-1}-1)}{\lambda _{\kappa }h(\theta )}\right) \text{ for all }%
\theta \in \lbrack 0,\bar{\theta}_{1}].
\]%
The equivalent statement is clear. \hfill $\square $

\subsection{Proof of Lemma 3}

\noindent It follows from our assumptions that the co-state variable takes
the following form:%
\begin{eqnarray*}
\text{For }0 &<&t<\theta _{1}:\pi =1\quad \mu _{t}=(\beta -\alpha )(\theta
_{1}-\theta _{2})e^{\lambda \lbrack g(t)-g(\theta _{1})]}; \\
\text{For }\theta _{1} &<&t<\theta _{2}:\pi =0\quad \mu _{t}=(\beta -\alpha
)(t-\theta _{2}); \\
\text{For }\theta _{2} &<&t<1:\pi =1\quad \mu _{t}=0.
\end{eqnarray*}%
So%
\begin{eqnarray*}
\text{For }0 &<&t<\theta _{1}:\pi =1\quad \gamma _{t}^{\ast }=\frac{%
e^{-\lambda g(t)}}{\kappa (1+\lambda h(t)(\theta _{2}-\theta _{1})e^{\lambda
\lbrack g(t)-g(\theta _{1})]}]}; \\
\text{For }\theta _{1} &<&t<\theta _{2}:\pi =0\quad \gamma _{t}^{\ast }=%
\frac{e^{-\lambda g(t)}}{\kappa (1+\lambda h(t)(\theta _{2}-t))}; \\
\text{For }\theta _{2} &<&t<1:\pi =1\quad \gamma _{t}^{\ast }=\frac{%
e^{-\lambda g(t)}}{\kappa }.
\end{eqnarray*}%
For $0<t<\theta _{1},$ if $\gamma _{t}>\gamma _{t}^{\ast },$ then%
\begin{eqnarray*}
e^{-\lambda g(t)} &>&\frac{e^{-\lambda g(t)}}{(1-\frac{\alpha }{\beta }%
)(1+\lambda h(t)(\theta _{2}-\theta _{1})e^{\lambda \lbrack g(t)-g(\theta
_{1})]}]}, \\
(1+\lambda h(t)(\theta _{2}-\theta _{1})e^{\lambda \lbrack g(t)-g(\theta
_{1})]}] &>&\kappa ^{-1}=1+\lambda h(\theta _{1})(\theta _{2}-\theta _{1}),
\\
h(t)e^{\lambda g(t)} &>&h(\theta _{1})e^{\lambda g(\theta _{1})}.
\end{eqnarray*}%
We observe that%
\[
\frac{d}{dt}h(t)e^{\lambda g(t)}=e^{\lambda g(t)}[h^{\prime }(t)+\lambda
h(t)^{2}].
\]%
As $h^{\prime }(t)<0,$ the derivative is negative iff
\[
h^{\prime }(t)+\lambda h(t)^{2}<0,
\]%
and, consequently, $h(t)e^{\lambda g(t)}$ needs to be decreasing on $%
[0,\theta _{1}]$ to achieve $h(t)e^{\lambda g(t)}>h(\theta _{1})e^{\lambda
g(\theta _{1})}.$ Recall that
\[
-h^{\prime }(t)=\sigma \exp \left( -[\sigma ^{2}(1-t)^{2}/8]\right) /\sqrt{%
2\pi }>0,
\]%
which is increasing with $t$ as likewise is $1/h(t).$ So the required
condition holds provided%
\[
\lambda <\eta (\sigma ):=-h^{\prime }(0)/h(0)^{2}=\sigma \exp \left( -\sigma
^{2}/8\right) /(\sqrt{2\pi }[2\Phi (\sigma /2)-1]^{2}),
\]%
yielding the claimed bound on $\lambda $. \hfill $\square $

\subsection{Proof of Proposition A}

Define $\kappa _{c}$ by
\[
\frac{\log \kappa _{c}^{-1}}{\kappa _{c}^{-1}-1}=\frac{1}{2}.
\]%
Take $\kappa \leq \kappa _{c}$ and put%
\[
\eta :=\frac{\kappa ^{-1}-1}{\lambda }\text{ and let }\bar{\theta}+\frac{%
\eta }{h(\bar{\theta})}=1.
\]%
The latter equation is equivalent to%
\[
\eta =(1-\bar{\theta})h(\bar{\theta}),
\]%
which is soluble since $h(0)>\eta .$ So $T(\bar{\theta})=1.$

Since $h=g^{\prime },$ $h^{\prime }<0$ and $h^{\prime \prime }<0$ on $[0,1],$
by Taylor's theorem applied to $\theta \leq \bar{\theta},$ there is some $%
\xi \in (\theta ,\theta +\eta /h(\theta ))$ with%
\[
g\left( \theta +\frac{\eta }{h(\theta )}\right) =g(\theta )+\frac{\eta }{%
h(\theta )}h(\theta )+\frac{\eta ^{2}}{2h(\theta )^{2}}h^{\prime }(\xi
_{\eta }).
\]%
Here $h^{\prime }(\xi _{\eta })>h^{\prime }(1)$ as $h^{\prime \prime }<0.$
Given the assumption about $\kappa ,$%
\[
\frac{\eta }{2}=\frac{1}{2}\frac{\kappa ^{-1}-1}{\lambda }>\frac{\log \kappa
^{-1}}{\kappa ^{-1}-1}\frac{\kappa ^{-1}-1}{\lambda }=\frac{\log \kappa ^{-1}%
}{\lambda }.
\]%
Given the assumptions about $\lambda $, since $h(\theta )\geq h(\bar{\theta}%
),$%
\[
\eta \frac{-h^{\prime }(1)}{h(\theta )^{2}}=\frac{\kappa ^{-1}-1}{\lambda }%
\frac{-h^{\prime }(1)}{h(\theta )^{2}}<\frac{\kappa ^{-1}-1}{\lambda }\frac{%
-h^{\prime }(1)}{h(\bar{\theta})^{2}}<1.
\]%
So%
\[
1+\frac{\eta }{h(\theta )^{2}}h^{\prime }(1)>0.
\]%
Combining gives%
\[
\frac{\eta }{2}+\frac{\eta }{2}\left( 1+\frac{\eta }{h(\theta )^{2}}%
h^{\prime }(1)\right) \geq \frac{1}{\lambda }\log \kappa ^{-1},
\]%
which implies the result via the Taylor expansion.\hfill $\square $

\subsection{Proof of Proposition B}

Here we employ the switching curve derived in the Hamiltonian formulation.
The equilibrating result follows from the equations asserting that switching
occurs when the state trajectory $\gamma _{t}$ crosses the switching curve $%
\gamma _{t}^{\ast }$, i.e.%
\[
\gamma _{\theta _{j}}=\gamma _{\theta _{j}}^{\ast }\text{ for }%
j=i-2,i-1,i,i+1.
\]

We first find the form of the costate variable $\mu _{t}.$ In view of the
alternations in policy we conclude as follows.%
\begin{eqnarray*}
\text{For }\theta _{i} &\leq &t\leq \theta _{i+1}:\pi =0\quad \mu
_{t}=(\beta -\alpha )\psi (t),\text{ for some linear }\psi (t):=(K-t). \\
\text{For }\theta _{i-1} &\leq &t\leq \theta _{i}:\pi =1\quad \mu
_{t}=(\beta -\alpha )\psi (\theta _{i})e^{\lambda \lbrack g(t)-g(\theta
_{i})]}. \\
\text{For }\theta _{i-2} &\leq &t\leq \theta _{i-1}:\pi =0\quad \mu
_{t}=(\beta -\alpha )(t-\theta _{i-1}+\psi (\theta _{i})e^{\lambda \lbrack
g(\theta _{i-1})-g(\theta _{i})]}).
\end{eqnarray*}%
W.l.o.g,, since $\gamma _{t}=\gamma _{t}^{\ast }$ at $t=\theta _{i-2},$ we
may rescale $\gamma _{t}$ so that $\gamma _{\theta _{i-2}}=1=\gamma _{\theta
_{i-2}}.$ So%
\begin{eqnarray*}
\noindent \text{1. for }\theta _{i} &\leq &t\leq \theta _{i+1}:\pi =0\quad
\gamma _{t}^{\ast }=\frac{e^{-\lambda g(t)}}{\kappa (1+\lambda h(t)\psi (t)]}%
, \\
\noindent \text{2. for }\theta _{i-1} &\leq &t\leq \theta _{i}:\pi =1\quad
\gamma _{t}^{\ast }=\frac{e^{-\lambda g(t)}}{\kappa (1+\lambda h(t)\psi
(\theta _{i})e^{\lambda \lbrack g(t)-g(\theta _{i})]}]}, \\
\noindent \text{3. for }\theta _{i-2} &\leq &t\leq \theta _{i-1}:\pi =0\quad
\gamma _{t}^{\ast }=\frac{e^{-\lambda g(t)}}{\kappa (1+\lambda h(t)(\theta
_{i-1}-t+\psi (\theta _{i})e^{\lambda \lbrack g(\theta _{i-1})-g(\theta
_{i})]})]}.
\end{eqnarray*}

The consistency of these at $\theta _{i-1}$ implies from (3) that%
\[
1=\gamma _{\theta _{i-1}}=\gamma _{\theta _{i-1}}^{\ast }=\frac{e^{-\lambda
g(\theta _{i-1})}}{\kappa (1+\lambda h(\theta _{i-1})\psi (\theta
_{i})e^{\lambda \lbrack g(\theta _{i-1})-g(\theta _{i})]}]}.
\]%
From $\gamma _{\theta _{i-1}}=1,$ it follows that
\[
\gamma _{t}=e^{-\lambda \lbrack g(t)-g(\theta _{i-1})]}\text{ for }\theta
_{i-1}\leq t\leq \theta _{i}\text{ (where }\pi =1).
\]%
Recalling that%
\[
\gamma _{\theta _{i}}=\gamma _{\theta _{i}}^{\ast },
\]%
equation (1) and the formula for $\gamma _{t}$ gives at $t=\theta _{i}$ that%
\[
\gamma _{\theta _{i}}=e^{-\lambda \lbrack g(\theta _{i})-g(\theta _{i-1})]}=%
\frac{e^{-\lambda g(\theta _{i})}}{\kappa (1+\lambda h(\theta _{i})\psi
(\theta _{i})]},
\]%
and so cross-multiplying,%
\[
\text{ }1=\frac{e^{-\lambda g(\theta _{i-1})}}{\kappa (1+\lambda h(\theta
_{i})\psi (\theta _{i})]}.
\]%
Now by equation (3)
\[
1=\gamma _{\theta _{i-1}}=\gamma _{\theta _{i-1}}^{\ast }=\frac{e^{-\lambda
g(\theta _{i-1})}}{\kappa (1+\lambda h(\theta _{i-1})\psi (\theta
_{i})e^{\lambda \lbrack g(\theta _{i-1})-g(\theta _{i})]}]}.
\]%
The last two results give%
\begin{eqnarray*}
\frac{e^{-\lambda g(\theta _{i-1})}}{\kappa (1+\lambda h(\theta _{i-1})\psi
(\theta _{i})e^{\lambda \lbrack g(\theta _{i-1})-g(\theta _{i})]}]} &=&\frac{%
e^{-\lambda g(\theta _{i-1})}}{\kappa (1+\lambda h(\theta _{i})\psi (\theta
_{i})]}, \\
h(\theta _{i-1})\psi (\theta _{i})e^{\lambda \lbrack g(\theta
_{i-1})-g(\theta _{i})]} &=&h(\theta _{i})\psi (\theta _{i}), \\
h(\theta _{i-1})e^{\lambda \lbrack g(\theta _{i-1})-g(\theta _{i})]}
&=&h(\theta _{i}),
\end{eqnarray*}%
and so finally%
\[
h(\theta _{i})e^{\lambda g(\theta _{i})}=h(\theta _{i-1})e^{\lambda g(\theta
_{i-1})},
\]%
as claimed.\hfill $\square $

\textbf{Symbols list:}

$\alpha ,\,\beta ,\kappa =1-\alpha /\beta :$ pay-for-performance
coefficients, \S 2 and 3.3;

$\pi :$ probability of sparing behaviour, \S 2;

$\gamma :$ valuation curves: $\gamma _{t}^{1}$ sparing throughout, $\gamma
^{\ast }$ switching curve, \S 3.3;

$\lambda :$ Poisson news-arrival intensity, \S 2;

$\theta :$ switching time, \S 2;

$\mu :$ co-state variable, \S 3.4;

$\varphi :$ integrating factor (associated with $\mu $), Appendix;

$h:$ discounting function (depending on time), \S 2 and 5;

$g(t):=\int_{0}^{t}h(u)$ $\mathrm{d}u;$

$X_{t}:$ economic state; $Y_{t}:$ its (noisy) observation, \S 2;

$\sigma :$ aggregate volatility, \S 3.2.

\end{document}